\documentclass{bmcart}

\usepackage{amsmath,amsthm,amssymb,amsopn,amsfonts,pdfpages,url, dsfont}
\usepackage[T1]{fontenc}
\usepackage[utf8]{inputenc}
\usepackage[linesnumbered,ruled,vlined]{algorithm2e}
\usepackage[colorlinks=true,linkcolor=cyan,citecolor=blue]{hyperref}
\usepackage{graphics,relsize}
\usepackage{graphicx}
\usepackage{multirow}
\usepackage{bbm}
\usepackage{authblk}
\usepackage{subcaption}
\usepackage{float}

\usepackage[utf8]{inputenc} 

\startlocaldefs
\endlocaldefs

\newcommand{\n}{\mathfrak{n}}
\newcommand{\SSS}{{\mathcal S}}
\newcommand{\N}{{\mathcal N}}

\newcommand{\st}{\mathfrak{str}}

\newcommand{\q}{\mathbbm{1}}
\newcommand{\dis}{\mathcal D}
\newcommand{\phh}{\hat{\varphi}}
\newcommand{\phS}{\hat{\varphi}_{\textup{SGM}}}
\newcommand{\qh}{\hat{q}}
\newcommand{\qhs}{\hat{q}_{\textup{SGM}}}

\newcommand{\dens}{\mathfrak{d}}
\newcommand{\Pry}{\mathbb{P}}
\newcommand{\e}{\mathbb{E}}

\newtheorem{theorem}{Theorem}

\newtheorem{definition}{Definition}

\begin{document}

\begin{frontmatter}

\begin{fmbox}
\dochead{Research}


\title{ The Phantom Alignment Strength Conjecture: Practical use of graph matching alignment strength to indicate a meaningful graph match  }


\author[
  addressref={aff1},                   
  corref={aff1},                       
  email={def@jhu.edu}   
]{\inits{D.E.}\fnm{Donniell E.} \snm{Fishkind}}
\author[
  addressref={aff2},
  email={fparker9@jhu.edu}
]{\inits{F.}\fnm{Felix} \snm{Parker}}
\author[
  addressref={aff1},
  email={hsawczu1@jhu.edu}
]{\inits{H.}\fnm{Hamilton} \snm{Sawczuk}}
\author[
  addressref={aff1},
  email={lmeng2@jhu.edu}
]{\inits{L.}\fnm{Lingyao} \snm{Meng}}
\author[
  addressref={aff3},
  email={ebridge2@jhu.edu}
]{\inits{E.}\fnm{Eric} \snm{Bridgeford}}
\author[
  addressref={aff1},
  email={dathrey1@jhu.edu}
]{\inits{A.}\fnm{Avanti} \snm{Athreya}}
\author[
  addressref={aff1},
  email={cep@jhu.edu}
]{\inits{C.E.}\fnm{Carey} \snm{Priebe}}
\author[
  addressref={aff4},
  email={vlyzinsk@umd.edu}
]{\inits{V.P.}\fnm{Vince} \snm{Lyzinski}}


\address[id=aff1]{
  \orgdiv{Department of Applied Mathematics and Statistics},             
  \orgname{Johns Hopkins University},          
  \city{Baltimore},                              
  \cny{MD}                                    
}
\address[id=aff2]{
  \orgdiv{Center for Systems Science and Engineering},             
  \orgname{Johns Hopkins University},          
  \city{Baltimore},                              
  \cny{MD}                                    
}
\address[id=aff3]{
  \orgdiv{Department of Biostatistics},             
  \orgname{Johns Hopkins University},          
  \city{Baltimore},                              
  \cny{MD}                                    
}
\address[id=aff4]{%
  \orgdiv{Department of Mathematics },
  \orgname{University of Maryland, College Park},
  \city{College Park},
  \cny{MD}
}



\end{fmbox}


\begin{abstractbox}

\begin{abstract} 
The alignment strength of a graph matching is a quantity that gives the practitioner a measure of the correlation of
the two graphs, and it can also give the practitioner a sense for whether the
graph matching algorithm found the true matching. Unfortunately, when a graph matching algorithm
fails to find the truth because of weak signal, there may be ``phantom
alignment strength'' from meaningless matchings that, by random noise,
have fewer disagreements than average (sometimes substantially fewer);
this  alignment strength may give the misleading appearance of significance.
A practitioner needs to know what level of alignment strength may be phantom alignment strength and what level
indicates that the graph matching algorithm obtained the true matching and is a meaningful measure of the
graph correlation. The {\it Phantom Alignment Strength Conjecture}
introduced here provides a principled and practical means to approach this issue. We provide empirical
evidence for the conjecture, and explore its consequences.
\end{abstract}




\end{abstractbox}
%

\end{frontmatter}



\section{Introduction \label{sec:intro}}

This paper is about graph matchability in practice. Specifically,
when given two graphs and an unobserved ``true'' bijection
(also called ``true matching'' or ``true alignment'') between their vertices, will
exact (i.e. optimal)
graph matching and approximate graph matching
algorithms provide us with the matching which is the
``truth''? How might we know in actual practice whether the ``truth'' has been found?
Our work is in response to
the latter question. The main contribution here is our
formulation of the
{\it Phantom Alignment Strength Conjecture}
in Section \ref{sec:ASC}, followed up in Section \ref{sec:ASC}
with the practical implications of this conjecture in deciding when alignment strength is high enough to indicate
truth. This conjecture is also interesting as a theoretical matter,
completely aside from its consequences.

Graphs (networks) are a commonly used data modality for encoding relationships, interactions, and dependencies in data in an incredibly broad range of the sciences and engineering; this includes sociology (e.g., social network analysis \cite{wasserman}),
neuroscience connectomics \cite{sporns_complex,vogelstein2019connectal},
biology (e.g., biological interaction networks \cite{szklarczyk2015string,tong2004global}), and automated
 knowledge discovery \cite{wang2017knowledge}, to name just a few application areas.

 The graph matching problem is, given two graphs with the same number of vertices, to find the bijection between the vertex sets that minimizes the number of adjacency
 ``disagreements'' between the graphs. Often there is an underlying ``true'' bijection that the graph matching is
 attempting to recover/approximate.
 Sometimes part of this true bijection is known a-priori, in which case minimizing the number of disagreements
 over the remainder of the bijection is called seeded graph matching. Graph matching and seeded graph matching are formally defined  in Section~\ref{sec:overview}.

 Graph matching and seeded graph matching are used in a wide variety of places, and we mention just a few.
 Information about the interactions amongst objects of interest is sometimes split across multiple networks or multiple layers of the same network \cite{kivela2014multilayer}.
In many applications, such as neuroscience connectomics where, for example, DT-MRI derived graphs can be generated by aligning scans to a common template before uncovering the underlying edge structure \cite{gray2012magnetic}, the vertices across networks or across layers are a priori aligned and identified.
These aligned vertex labels can then be used to create joint network inference procedures that can leverage the signal across multiple networks for more powerful statistical inference \cite{levin2017central, chen2016joint, arroyo2019inference, durante2018bayesian}.
In many other applications, the vertex labels across networks or across layers are unknown or noisily observed.
Social networks provide a canonical example of this, where common users across different social network platforms may use different user names and their user profiles may not be linked across networks.
Discovering this latent correspondence (in the social network example, this is anchoring profiles to a common user across networks) is a key inference task \cite{lin2010layered,yartseva2013performance} for leveraging the information across networks for subsequent inference, and it is a key consideration for understanding the degree of user anonymity \cite{ding2010brief} across platforms.

 For a thorough survey of the relevant graph matching literature, see \cite{ConteReview,
 foggia2014graph,Emmert-Streib2016-st}.

The graph matching problem is computationally complex. Indeed, the simpler graph isomorphism problem has been shown to be of quasi-polynomial complexity \cite{babai2016graph}. Allowing  loopy, weighted, directed graphs makes graph matching equivalent to the NP-hard quadratic assignment problem. Due to its practical importance and computational difficulty, a large branch of the graph matching literature is devoted to developing algorithms to efficiently, but approximately, solve the graph matching problem; see, for example, \cite{FAP,umeyama,singh2007pairwise,zaslavskiy2009path,zhou,FAQ,zhang2016final,feizi2016spectral,heimann2018regal} among myriad others.

Somewhat dual to the algorithmic development literature, a large branch of the modern graph matching literature is devoted to theoretically exploring the question of graph matchability, also called graph de-anonymization; this is the
question of determining when there is enough signal present for graph matching to recover the ``true'' bijection.
Many of the recent papers in this area have introduced latent alignment across graphs by correlating the edges across networks between common pairs of vertices, focusing on understanding the phase transition between matchable and non-matchable networks in terms of the level of correlation across networks and/or the sparsity level of the networks; see, for example,
\cite{pedarsani2011privacy,lyzinski2014seeded,cullina2016improved,rel,cullina2017exact,sussman2018matched,cullina2019partial,fan2019spectral,ding2020efficient,mossel2020seeded}.

In \cite{fishkind2019alignment}, a novel measure of graph
correlation between two random graphs
called {\it total correlation} is introduced; it
is neatly partitioned into an inter-graph contribution (the ``edge correlation'' that had been the previous focus in the literature) and a novel intra-graph contribution. Furthermore, they introduce
a statistic called {\it alignment strength}, which
is $1$ minus a normalized count of the number of
disagreements in an optimal/true graph match; they prove
under mild conditions that
alignment strength is a strongly consistent estimator
of total correlation.
Experimental results in \cite{fishkind2019alignment} suggest that the matchability phase transition, as well as the complexity of the problem, is a function of this more nuanced total correlation rather than simply the cross-graph edge correlation/edge sparsity that had been the previous focus in the literature.

Analyses mining the matchability phase transition in the literature that also have considered similarity across generative network models beyond simple sparsity have thus far focused on simple community-structured network models
\cite{onaran2016optimal,shirani2018matching,lyzinski2016information}, or have proceeded by removing the heterogeneous within-graph model information and simply using the across graph edge correlation \cite{lyzinski2017matchability}.
Recently, there have been numerous papers in the literature at the interface between algorithm development and mining matchability phase-transitions; see, for instance, \cite{barak2019nearly,mossel2020seeded,ding2020efficient}.
A common theme of many of these results is that, under assumptions on the across graph edge-correlation and network sparsity, algorithms are designed to efficiently (or approximately efficiently) match graphs with corresponding theoretical guarantees on the
performance of the algorithms in recovering the latent alignment.

However, the question remains how a practitioner
knows in practice whether or not a graph matching has successfully recovered the truth. This issue is not resolved by asymptotic
analysis with hidden constants. Nor, in general, are the underlying parameters known to the practitioner.
It seems that the graph alignment statistic is a very
natural metric to use in
deciding if the truth is found. Unfortunately, when there is an absence of signal, an optimal (or approximately optimal) graph matching will find spurious and random
alignment strength due to chance. Indeed, this meaningless
alignment strength can be high and misleading. How do we
gauge whether or not  it is high enough to signal that truth is found?

After formally defining seeded graph matching and alignment
strength in Section \ref{sec:overview} and defining the
correlated Bernoulli random graph model (and attendant parameters) in Section \ref{sec:cbrg},
we then address this issue with our Phantom Alignment
Strength Conjecture in Section \ref{sec:ASC}, and in the
ensuing discussion in Section \ref{sec:ASC}. Then, in Section \ref{sec:evidence}, we
present empirical evidence for the conjecture
using synthetic and real data, and comparing to theoretical results; Section \ref{sec:evidence} begins with a thorough summary.
This is followed in Section \ref{sec:notable}
by notable mentions, and future directions.

\section{Seeded graph matching, alignment strength \label{sec:overview}}

\indent \indent In the seeded graph matching setting, we are given two simple graphs,
say they are $G_1=(V_1,E_1)$ and $G_2=(V_2,E_2)$, such that $|V_1|=|V_2|$, denote the number
of vertices $\n:=|V_1|$.
Let $\varPi$ denote the set of all
bijections $V_1\rightarrow V_2$.
It is usually understood that there
exits a ``true'' bijection $\varphi^* \in \varPi$ which represents a natural correspondence between the vertices in $V_1$
and the vertices~in~$V_2$; for example, $V_1$ and $V_2$ might be the same people, with $E_1$ indicating which pairs exchanged
emails and $E_2$ indicating pairs that communicated in a different medium.
Or $G_1$ may be the electrical connectome (brain graph) of a worm and
$G_2$ might be the chemical connectome of the same worm, both graphs sharing the same vertex set of neurons.
The vertex set $V_1$ is partitioned into two disjoint sets, $\SSS$ ``seeds'' (possibly empty) and $\N$ ``nonseeds,'' denote
$s:=|\SSS|$ and $n:=|\N|$. (When $s=0$ this is the conventional graph matching problem.)
The graphs $G_1$ and $G_2$ are observed, and the values of $\varphi^*$ are observed on the set of seeds $\SSS$,
however the values of  $\varphi^*$ are not observed on the nonseeds $\N$, and one of several important tasks is to estimate $\varphi^*$.

 Let $\varPi^S$ denote the set of all
bijections $V_1\rightarrow V_2$ that agree with $\varphi^*$ on the seeds $\SSS$. For any $\varphi \in \varPi^S$, its
{\it match ratio} is defined to be  $\frac{1}{n}|\{ v \in \N : \varphi(v)=\varphi^*(v) \}|$, i.e.
the fraction of the nonseeds that are correctly matched by $\varphi$. (It is common to multiply the
match ratio by $100$ to express it as a percentage.)

For any set~$V$, let ${V \choose 2}$ denote the set of two-element subsets of $V$; for each $i=1,2$ and any $\{u,v\} \in {V_i \choose 2}$
let $u \sim_{G_i} v$ and  $u \not \sim_{G_i} v$ denote adjacency and, respectively, nonadjacency of $u$ and $v$ in $G_i$.
Next, let~$\q$~denote the indicator function for its subscript.
Given any $\varphi \in \varPi$, we define the {\it full number of disagreements through} $\varphi$ to be
\begin{eqnarray} \label{eqn:dis}
\dis' (\varphi):= \sum_{ \{u,v\} \in {V_1 \choose 2}  } \left  (  \q_{[u \sim_{G_1} v] \wedge [\varphi(u) \not \sim_{G_2} \varphi(v)]}+
 \q_{[u \not \sim_{G_1} v] \wedge [\varphi(u) \sim_{G_2} \varphi(v)]} \right )
\end{eqnarray}
and, given any $\varphi \in \varPi^S$, we define  the {\it restricted number of disagreements through} $\varphi$ to be
\begin{eqnarray} \label{eqn:dis2}
\dis (\varphi):= \sum_{ \{u,v\} \in {\N \choose 2}  } \left  (  \q_{[u \sim_{G_1} v] \wedge [\varphi(u) \not \sim_{G_2} \varphi(v)]}+
 \q_{[u \not \sim_{G_1} v] \wedge [\varphi(u) \sim_{G_2} \varphi(v)]} \right ).
\end{eqnarray}
The seeded graph matching problem is to find
\begin{eqnarray} \phh \in \arg \min _{\varphi \in \varPi^S} \dis'(\varphi), \label{eqn:optim}
\end{eqnarray}
and the idea is that $\phh$ is an estimate for the true bijection $\varphi^*$. Unfortunately, except in the smallest instances,
computing $\phh$ is intractable. A state-of-the-art algorithm SGM from \cite{FAP} is commonly used to approximately solve the optimization
problem in (\ref{eqn:optim}), and we denote its output $\phS \  (\in \varPi^S)$, and it is an approximation of $\phh$ and, hence,
an approximation of $\varphi^*$.
For any $\varphi \in \varPi^S$,
the {\it full alignment strength} $\st' (\varphi)$ and the {\it restricted alignment strength} $\st (\varphi)$ are defined as
\begin{eqnarray}
\st' (\varphi):= 1 -\frac{\dis' (\varphi) }{\frac{1}{\n !}\sum_{\phi \in \varPi} \dis' (\phi )  } \ \ \ \mbox{  and  }  \ \ \
\st (\varphi):= 1 -\frac{\dis (\varphi) }{\frac{1}{n !}\sum_{\phi \in \varPi^S} \dis (\phi )  } .
\label{eqn:as}
\end{eqnarray}

Although the denominators of (\ref{eqn:as}) have exponentially many summands, alignment strength is easily computed
as follows. For $i=1,2$, define the full density of $G_i$ as $\dens' G_i:=\frac{|E_i|}{{\n \choose 2}}$ and
the restricted density of $G_i$ as $\dens G_i=$ the number of edges of $G_i$ induced by $\N$, divided by ${n \choose 2}$.
It holds that
\begin{eqnarray}
\st' (\varphi) & = & 1 -\frac{  \dis' (\varphi)/{\n \choose 2} }{\dens' G_1 (1-\dens' G_2)+ (1-\dens' G_1) \dens' G_2 } \ \ \  \nonumber \mbox{ and } \\
\st (\varphi)  & = & 1 -\frac{  \dis (\varphi)/{n \choose 2} }{\dens G_1 (1-\dens G_2)+ (1-\dens G_1) \dens G_2 };
\label{eqn:asf}
\end{eqnarray}
see \cite{fishkind2019alignment} for the derivation of (\ref{eqn:asf}) from (\ref{eqn:as}).

The importance of alignment strength to a practitioner is twofold:

First, the alignment strength of $\varphi^*$ (and its proxies $\phh$ and $\phS$) may be thought of as a measure of how
similar the structure of the graphs $G_1$ and $G_2$ are through the ``true'' bijection; indeed,
if the number of disagreements under $\varphi^*$  (and its proxies $\phh$ and $\phS$) is about equal to the average over all bijections
then its alignment strength is near $0$ (as clearly seen from the definition in (\ref{eqn:as})) and, at the other extreme,
if $\varphi^*$ (and its proxies $\phh$ and $\phS$) is nearly an isomorphism between $G_1$ and $G_2$ then its alignment strength is near~$1$.
It was proven in \cite{fishkind2019alignment} that the full alignment strength of the ``true'' bijection
$\st' (\varphi^*)$ is a strongly consistent estimator of $\varrho_T$, which is a parameter called the
{\it total correlation} between the two graphs $G_1$ and $G_2$, defined in Section~\ref{sec:cbrg}.

Another way that alignment strength is of much importance to a practitioner is in
providing confidence that $\phS$ or $\phh$ is a good estimate of $\varphi^*$, the ``truth.''
If $\st (\phS)$ or $\st (\phh)$ is high enough then we may be confident that a meaningful match capturing similar graph
structure has been found, and therefore $\phS$ or $\phh$ is approximately or exactly $\varphi^*$.
But, how high is high enough?

Indeed, these issues in the use of alignment strength become vastly more complicated by
the possibility of {\it phantom alignment strength}.
This is a phenomenon that occurs when, in the presence of weak signal,  meaningless matchings have many fewer disagreements
than average (sometimes very substantially fewer) due to random noise,
and $\phh$ and/or $\phS$ is one of these meaningless matchings---optimal in the optimization problem, but meaningless
as estimates of $\varphi^*$. Indeed, the alignment strength of $\phh$ and/or $\phS$
may be elevated enough to give the misleading appearance of significance when, in reality,
they don't at all resemble $\varphi^*$. This will be illustrated in Section \ref{sec:evidence}.

The purpose of this paper is to give a principled, practical
means of approaching the decision of what level of alignment strength for $\phh$ and/or $\phS$
indicates that they are a good approximation of $\varphi^*$, in which case the alignment strength
reflects the amount of meaningful similar structure
between $G_1$ and $G_2$ ---beyond the random similarity  between completely unrelated graphs.\\

(A note on terminology:
We define both {\it full} alignment strength and {\it restricted} alignment strength since each will end up being important at a different time.
The Phantom Alignment Strength Conjecture of Section \ref{sec:ASC} requires restricted alignment strength specifically; indeed, since
full alignment strength includes the seeds, this would dilute the desired effect, falsifying the conjecture conclusion.
However, after we have confidence that our graph matching is the true matching, it is then full
alignment strength that will be a better estimator of total correlation introduced in Section \ref{sec:cbrg}.)

\section{The correlated Bernoulli random graph model \label{sec:cbrg}}

\begin{definition}
Given positive integer $\n$, vertex set $V$ such that $|V|=\n$, the
parameters~of~the~\textup{correlated Bernoulli random graph model} are Bernoulli parameters $p_{\{u,v\}} \in [0,1]$ for each
$\{u,v\} \in {V \choose 2}$, and an edge correlation parameter $\varrho_e \in [0,1]$.
The pair of random graphs $(G_1,G_2)$ have a \textup{correlated Bernoulli random graph distribution} when as follows:
$G_1$ and $G_2$ each have vertex set $V$.
For each $\{u,v\} \in {V \choose 2}$, and each $i=1,2$, the probability of $u \sim_{G_i} v$ is the Bernoulli parameter
$p_{\{u,v\}}$, and the Pearson correlation for random variables $\q_{v \sim_{G_1} w}$ and $\q_{v \sim_{G_2} w}$
is~the~edge~correlation parameter $\varrho_e$. Other than these dependencies, the rest of the adjacencies are independent.
\end{definition}
The distribution of the pair of random graphs $G_1,G_2$ is determined by the above (see \cite{fishkind2019alignment}). Of~course, the identity
function  is the ``true'' matching $\varphi^*$ between $G_1$ and $G_2$.

(If the Bernoulli parameters are all equal, then the random graphs $G_1$ and $G_2$ are each said to be Erdos-Renyi, so the correlated Erdos-Renyi random graph model is a special case of the correlated Bernoulli random graph model.)

Important functions of the model parameters are as follows. The Bernoulli mean and Bernoulli variance are, respectively, defined as
\begin{eqnarray*}
\mu:= \frac{ \sum_{\{u,v\}\in {V \choose 2}}p_{\{u,v\}}}{{\n \choose 2}},
\ \ \ \ \sigma^2:=\frac{\sum_{\{u,v\}\in {V \choose 2}}(p_{\{u,v\}}-\mu)^2}{{\n \choose 2}} .
\end{eqnarray*}
Assume that $\mu$ is not equal to $0$ nor $1$.
The {\it heterogeneity correlation} is defined in
\cite{fishkind2019alignment} as
\begin{eqnarray} \varrho_h := \frac{\sigma^2}{\mu (1-\mu)}; \label{eqn:hc}
\end{eqnarray}
it is in the unit interval $[0,1]$; see \cite{fishkind2019alignment}. Also pointed out in \cite{fishkind2019alignment} is that $\varrho_h$ is $0$ if and only if all Bernoulli parameters are equal (i.e. the graphs are Erdos-Renyi) and $\varrho_h$ is $1$ if and only if all Bernoulli parameters are $\{0,1\}$-valued. In particular, if $\varrho_h$ is~$1$ then
$G_1$ and $G_2$ are almost surely isomorphic.
The {\it total correlation} $\varrho_T$ is defined in
\cite{fishkind2019alignment}
to satisfy the relationship
\begin{eqnarray} (1-\varrho_T)=(1-\varrho_h)(1-\varrho_e) . \label{eqn:tc}
\end{eqnarray}

In the following key result, Theorem \ref{thm:theresult}, which was
 proved in \cite{fishkind2019alignment}, let us consider a probability space that
incorporates correlated Bernoulli random graph distributions for each of
the number of vertices $\n=1,2,3,\ldots$. Thus, the parameters are
functions of $\n$, but to prevent notation clutter we omit notating the dependence on $\n$. The symbol
$\stackrel{a.s.}{\rightarrow}$ denotes almost sure convergence.

\begin{theorem} \label{thm:theresult} Suppose $\mu$ is bounded away from $0$ and $1$, over all $\n$. Then it holds that
$\st' (\varphi^*) - \varrho_T \stackrel{a.s.}{\rightarrow} 0$.
\end{theorem}

Theorem \ref{thm:theresult} together with
Equation \ref{eqn:tc} shows that the alignment strength
of the true bijection captures (asymptotically) an
underlying correlation between the random graphs
that can be neatly (and symmetrically,
per Equation \ref{eqn:tc}) partitioned
into a inter-graph contribution
(edge correlation) and an intra-graph contribution
(heterogeneity correlation).\\

Next, instead of considering a sequence of correlated Bernoulli random graphs, let us dig down deeper
one probabilistic
level. Specifically, suppose that for
each $\{u,v\} \in {V \choose 2}$ there exists an
interval-$[0,1]$-valued distribution $F_{\{u,v\}}$
such that the Bernoulli parameter  $p_{\{u,v\}}$
(in the correlated Bernoulli random graph model)
is an independent random variable with distribution
$F_{\{u,v\}}$. Denote the mean of this distribution
$\mu_{F_{\{u,v\}}}$, denote
the variance of this distribution
$\sigma^2_{F_{\{u,v\}}}$, and (if we have
$\mu_{F_{\{u,v\}}}$ not $0$ nor $1$)
define the heterogeneity
correlation of the distribution to be
\begin{eqnarray}
\varrho_{F_{\{u,v\}}}:=\frac{\sigma^2_{F_{\{u,v\}}}}{\mu_{F_{\{u,v\}}}(1-\mu_{F_{\{u,v\}}})} \ \ .
\end{eqnarray}

\begin{theorem} \label{thm:new}
Given an edge correlation parameter $\varrho_e \in [0,1]$ and, for each $\{u,v\}\in {V \choose 2}$, given a
$[0,1]$-valued distribution $F_{\{u,v \}}$ such that the
Bernoulli parameter $p_{\{u,v\}}$ is independently distributed as $F_{\{u,v \}}$, then the
distribution of the associated correlated Bernoulli
random graphs $(G_1,G_2)$ is completely specified by  $\varrho_e$ and,  for all  $\{u,v\}\in {V \choose 2}$,
the values of
$\mu_{F_{\{u,v\}}}$ and $\varrho_{F_{\{u,v\}}}$.
\end{theorem}
\noindent {\bf Proof:} Consider any
$\{u,v\}\in {V \choose 2}$;
the Bernoulli coefficient $p_{\{u,v\}}$, call it $X$, has distribution~$F_{\{u,v \}}$.
For any $p \in [0,1]$, conditioning on $X=p$,
the joint probabilities of combinations of
$u,v$ adjacency in $G_1,G_2$ are
computed in a straightforward way (see \cite{fishkind2019alignment} Appendix A) in the table:
\begin{eqnarray}
\begin{array}{c||cc}
         &    u\sim_{G_2}v          &  u\not \sim_{G_2}v         \\ \hline \hline
u\sim_{G_1}v     &  p^2+\varrho_e p(1-p)     &  (1-\varrho_e)p(1-p)  \\
u \not \sim_{G_1}v     &  (1-\varrho_e)p(1-p)      & (1-p)^2+\varrho_e p(1-p)
\end{array}
\end{eqnarray}
Probabilities of these adjacency
combinations, relative to the underlying
distribution $F_{\{u,v \}}$,
are computed by integrating/summing the conditional
probabilities (in table) times the density/mass of
$F_{\{u,v \}}$, obtaining
\begin{eqnarray*}
& & \Pry [u \sim_{G_1}v \mbox{ and } u \not \sim_{G_2}v ]=\Pry [u \not \sim_{G_1}v \mbox{ and } u  \sim_{G_2}v ] \\
& = & (1-\varrho_e)( \e X - \e X^2) \\
& = & (1-\varrho_e)(\e X -(\e X)^2 -\e X^2 + (\e X)^2) \\
& = & (1-\varrho_e)[ \mu_{F_{\{u,v\}}} (1-\mu_{F_{\{u,v\}}} ) \  - \ \sigma^2_{F_{\{u,v\}}} ] \\
& = &  \mu_{F_{\{u,v\}}} (1-\mu_{F_{\{u,v\}}} ) (1-\varrho_e)
(1-\varrho_{F_{\{u,v\}}}).
\end{eqnarray*}
Then, for each $i=1,2$,
because $\Pry [u \sim_{G_i}v]= \e X =\mu_{F_{\{u,v\}}} $
we have all four adjacency combinations as functions of
$\mu_{F_{\{u,v\}}}$ and $\varrho_{F_{\{u,v\}}}$. The result
follows from the independence across all pairs of vertices.
$\qed$

In the Phantom Alignment Strength Conjecture we assume
all distributions  $F_{\{u,v\}}$ are the same, call the common distribution $F$. Note that Bernoulli mean $\mu$
and heterogeneity correlation $\varrho_h$ are
now random variables, and if $\n$ is large,
then $\mu$ and $\varrho_h$ will respectively be good estimators of $\mu_F$ and $\varrho_F$. {\bf
A very important
consequence of Theorem \ref{thm:new} is that
the only information that matters
regarding $F$ is contained (well-estimated)
in the quantities
$\mu$ and $\varrho_h$. }

\section{Phantom Alignment Strength Conjecture, consequences  \label{sec:ASC}}

In this section, we propose the Phantom Alignment Strength Conjecture, which is the central purpose
of this paper. We then discuss its consequences; the conjecture gives us a principled and practical way to
decide if we should be convinced that the output of a graph matching algorithm
well-approximates the true matching.

Henceforth we use the term {\it alignment strength} to refer to the restricted alignment strength.

Consider correlated Bernoulli random graphs $G_1,G_2$ such that there are a
``moderate'' number $n$ of nonseed vertices (say $ n \geq 300$),
$s$ seeds (selected discrete uniformly from the $\n:=n+s$ vertices), and Bernoulli parameters are independently
realized from any fixed $[0,1]$-valued distribution with moderate mean $\mu'$ (say $.05 < \mu' < .95$).
The {\it Phantom Alignment Strength Conjecture} states
that, subject to caveats, as discussed in Section \ref{sec:caveat}, there exists a {\it phantom alignment strength value}
$\qh \equiv \qh (n,s,\mu') \in [0,1]$
such that $\st (\phh)$
has ``negligible'' variance and is approximately a function of the total correlation $\varrho_T$ and, specifically, it holds that,
with ``high probability,''
\begin{eqnarray} \label{eqn:one}
\st (\phh) \approx \left \{  \begin{array}{rl}  \varrho_T & \mbox{ if } \varrho_T > \qh; \mbox{ in which
case } \phh = \varphi^*  \\ \qh & \mbox{ if } \varrho_T \leq \qh; \mbox{ in which
case } \phh \mbox{ is ``very different'' from } \varphi^* \end{array}
 \right . .
\end{eqnarray}
Moreover, the conjecture states that, when using the seeded graph matching algorithm SGM of~\cite{FAP},
(given $n,s, \mu'$, as above) then there exists
$\qhs \equiv \qhs (n,s,\mu') \in [0,1]$ such that $\qhs \geq \qh$, and
$\st (\phS)$ has ``negligible'' variance and is approximately a function of the total correlation $\varrho_T$ and,
specifically, it holds that, with ``high probability,''
\begin{eqnarray} \label{eqn:two}
\st (\phS) \approx \left \{  \begin{array}{rl}  \varrho_T & \mbox{ if } \varrho_T > \qhs; \mbox{ in which
case } \phS = \varphi^*  \\ \qh & \mbox{ if } \varrho_T \leq \qhs;  \mbox{ in which
case } \phS \mbox{ is ``very different'' from } \varphi^*  \end{array}
 \right . .
\end{eqnarray}

Note that both $\st (\phh)$ and $\st (\phS)$ are conjectured to be an
approximately piecewise linear function of $\varrho_T$; two pieces, one piece with slope $0$ and one piece with slope $1$.
However, $\st (\phh)$ is continuous and shaped like a hockey stick (see Figure \ref{fig:52a:g}), whereas for
$\st (\phS)$ there can be a discontinuity (see Figure \ref{fig:52a:b}); but the function value of the linear portion with slope $0$ is
the same for $\st (\phS)$ as it is for $\st (\phh)$, namely it is the   phantom alignment strength value
$\qh$.

There are important consequences of the Phantom Alignment Strength Conjecture for the practitioner.
Suppose that a practitioner has two particular graphs $G_1,G_2$ with $n$ nonseed vertices and $s$ seeds
that can be considered as realized from a correlated Bernoulli random graph model,
and the practitioner wants to seeded graph match them, computing $\phS$ as
an approximation of the true matching $\varphi^*$.
How can the practitioner tell if $\phS$ is $\varphi^*$?
This conjecture provides a principled, practical mechanism.
The practitioner should realize two independent Erdos-Renyi graphs $H_1$ and $H_2$
with $n$ nonseed vertices, $s$~seeds, and adjacency probability parameter $p$
equal to the combined density of $G_1$ and $G_2$.
Then use SGM to seeded graph match $H_1$ and $H_2$,
and the alignment strength of the bijection (between $H_1$ and $H_2$) is approximately $\qh \equiv \qh (n,s,\mu)$, since the
total correlation in generating $H_1$ and $H_2$ is $0$, by design.
Then, when subsequently seeded graph matching $G_1$ and $G_2$, if $\st (\phS)$ is greater than some predetermined and fixed $\epsilon>0$
above~$\qh$, then that would indicate that $\phS = \varphi^*$ and, if
$\st (\phS)$ is less than this, then there is no confidence that $\phS$ is $\varphi^*$.
Moreover, in the former case the practitioner can have confidence in approximating $\st (\phS) \approx \varrho_T$, and
in the latter case there wouldn't be confidence in this approximation. (In the
former case, note that the full alignment strength $\st' (\phS)$ would then be an even better estimate
of $\varrho_T$.)

(If some of the model assumptions are violated and the
Bernoulli mean of $G_1$ may be different from $G_2$, then
it may be better not to combine their densities,
but rather to realize $H_1$ and $H_2$
as Erdos-Renyi graphs with respective adjacency
parameter equal to their respective densities.)

\section{Empirical evidence in favor of the Phantom Alignment Strength Conjecture \label{sec:evidence}}

In this section we provide empirical evidence for the Phantom Alignment Strength Conjecture.

 A summary
is as follows:

We begin in Section
\ref{sec:small} with a scale small enough
($n$ is just on the order of tens) to solve seeded
graph matching and attain
optimality.
Although the Phantom Alignment Strength
Conjecture does not apply because $n$ is so small,
we nonetheless see many ingredients of the conjecture. Then,
in Section~\ref{sec:broken}, we use synthetic data on a scale for the conjecture to be applicable, and we empirically demonstrate the conjecture for many types of Bernoulli parameter distributions; unimodal, bimodal, symmetric, skewed, etc. The SGM algorithm is employed for seeded graph matching, since exact optimality is unattainable in practice.

In Section \ref{sec:comp}, the alignment
strength of completely uncorrelated Erdos-Renyi graphs
(graph matched with SGM, using no seeds), taken
as a function of $n$, is empirically demonstrated to be
the same order of growth (in terms of $n$) as the theoretical bound for matchability (as a function of $n$), which suggests that the two quantities are the same, in excellent accordance
with the conjecture.

Then, in Section \ref{sec:block}, we observe
that when there is
block structure and differing distributions for the Bernoulli
parameters by block (thus the conjecture
hypotheses are not adhered to) then the conjecture's claims
may fail to hold, to some degree.
Nonetheless, there is still a phantom
alignment strength that allows for a procedure similar to what we recommend
in Section \ref{sec:ASC} to be successfully used
for deciding when alignment strength is significant enough to indicate that the seeded graph matching has found the truth.

Real data is then used for demonstration in Section \ref{sec:noisy} and Section \ref{sec:modalities}.

Specifically, in Section
\ref{sec:noisy}, we use a human connectome at many different
resolution levels, and graph match it to a manually noised
copy of itself.

Then, in Section \ref{sec:modalities}, we consider
several pairs of real-data graphs (titled Wikipdeia, Enron, and C Elegans) whose vertices are the same objects, and the
adjacencies in each pair of graphs
represent relationships between the objects across two different modalities.

All of these experiments serve as strong empirical evidence for the Phantom Alignment Strength Conjecture,
and motivate its use.

\subsection{Of hockey sticks and phantom alignment strength \label{sec:small} }

We begin with an experiment in which the value of $n$ is well below what is required in the statement of the
Phantom Alignment Strength Conjecture. However, $n$ is small enough here to enable us to compute $\phh$ exactly, using the integer programming formulation from~\cite{fishkind2019alignment}. We will be able to see many
features of the Phantom Alignment Strength Conjecture, and we will also see that phantom alignment
strength is not just an artifact of the SGM algorithm.

\begin{figure}[h!]
	\centering
	\includegraphics[width=4.5in]{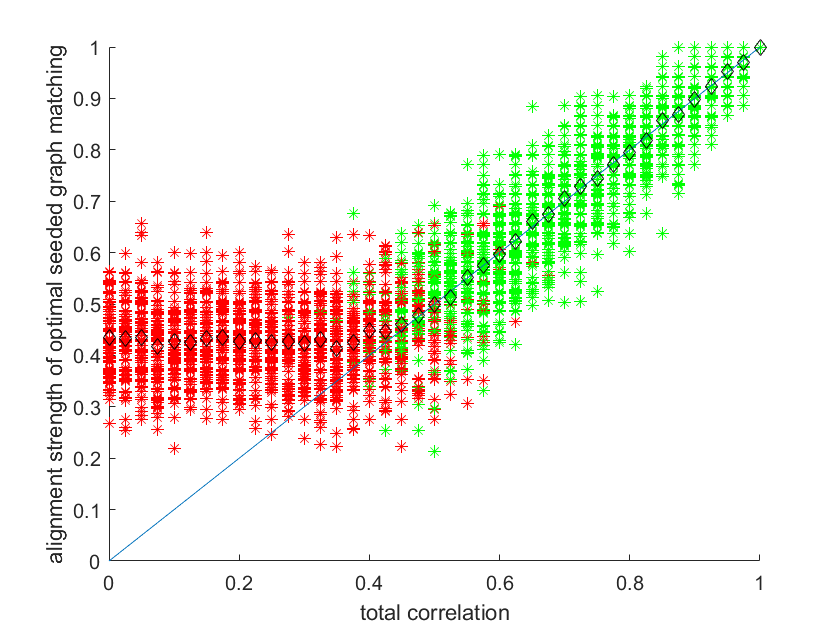}
    \caption{For each $\varrho_e$ from $0$ to $1$ in increments of $.025$, alignment strength of $\phh$
    for $100$ independent realizations when all Bernoulli probabilities were $0.5$ (in particular, $\varrho_T=\varrho_e$),
    with  $n=15$ nonseeds, $s=15$ seeds, a green asterisk if $\phh=\varphi^*$, else a red asterisk.}		
	\label{fig:exper1}
\end{figure}

For each value of $\varrho_e$ from $0$ to $1$ in increments of $.025$, we did $100$
independent repetitions of the following experiment. We realized a pair of correlated
Bernoulli random graphs on $\n=30$ vertices with edge correlation $\varrho_e$ and, for
each pair of vertices, the associated Bernoulli parameter~was~$0.5$.
(In particular, the graphs are correlated Erdos-Renyi.)
Since here $\sigma^2=0$,
we have that $\varrho_h=0$, and thus $\varrho_T=\varrho_e$.
We discrete uniform randomly chose $s=15$ seeds, so there were $n=15$ nonseeds.
For each experiment, we solved the seeded graph matching problem to optimality (indeed, $n=15$ is small
enough to do so), obtaining $\phh$. If it happened that $\phh=\varphi^*$ then we plotted a green
asterisk in Figure \ref{fig:exper1} for the resulting alignment strength $\st (\phh)$ against
the total correlation $\varrho_T$ and, if $\phh \ne \varphi^*$, we plotted a red asterisk
for the resulting alignment strength $\st (\phh)$ against the total correlation $\varrho_T$. The
black diamonds in Figure~\ref{fig:exper1} are the mean alignment strengths for the
$100$ repetitions, plotted for each value of $\varrho_e$.

It is readily seen from Figure \ref{fig:exper1}
that the variance for the
alignment strength of $\phh$ is quite high, which is reason
to not formulate the Phantom Alignment Strength Conjecture until $n$ is much larger. Other that this, observe that
if we substitute
``{\bf mean} of the alignment strength of $\phh$''
into the conjecture in place of ``alignment strength
 of $\phh$" then the conjecture would hold here.
Indeed, when $\varrho_T > \ \approx 0.44 \equiv \qh$ we very generally had that  $\phh=\varphi^*$, and
when $\varrho_T \leq  \ \approx 0.44$ we very generally had that  $\phh \ne \varphi^*$. (This boundary is not sharp, but is close.)
Also, note that when $\varrho_T >  \ \approx 0.44$, the mean of the alignment strength was approximately equal to $\varrho_T$.
Furthermore, when $\varrho_T \leq  \ \approx 0.44$, we see that
 the (mean) alignment strength of $\phh$ is the  phantom alignment strength (mean) of $ \ \approx 0.44$. Indeed, in this latter case, the alignment strength of $\phh$ is a misleading high value, and is not meaningful.

\subsection{Of hockey sticks and broken hockey sticks \label{sec:broken}}

In this section, we use synthetic data that meets the hypotheses of the Phantom Alignment Strength Conjecture.
Our setup was as follows.
We chose the number of nonseeds to be $n=1000$, and we repeated an
experiment for {\bf all combinations} of the following:
\begin{itemize}
\item Each pair of Beta distribution parameters $\alpha,\beta$ listed in the following table:
\[
\begin{array}{c|cc}
           & \alpha & \beta  \\ \hline
    \mbox{Pair A} &   1  &  1  \\
    \mbox{Pair B} &   0.5  &  0.5 \\
    \mbox{Pair C} &   2  &  2  \\
    \mbox{Pair D} &   5  &  1  \\
    \mbox{Pair E} &   2  &  5
\end{array}
\]
\item Each  $\mu'=$(mean of the scaled/translated Beta distribution) from $.1$ to $.9$ in increments of $.1$,
\item Each number of seeds $s= 0,10,20,50,250,1000$,
\item Each value of edge correlation $\varrho_e$ from $0$ to $1$ in increments of $0.025$,
\item Each value of $\delta$ from $0$ to $\delta_{max}:=\min \{ \frac{\alpha+\beta}{\alpha}\mu',\frac{\alpha+\beta}{\beta}(1-\mu') \}$
in increments of $\frac{1}{10}\delta_{max}$.
\end{itemize}
For each combination of the above, we realized a pair of correlated Bernoulli random graphs on $n+s$ vertices, with edge correlation~$\varrho_e$
and, for each pair of vertices, the associated Bernoulli parameter  was independently
realized from the distribution $\delta \cdot \textup{Beta} (\alpha,\beta)+ \mu'-\delta\frac{\alpha}{\alpha+\beta}$.
Note that
\begin{itemize}
\item The distribution $\delta \cdot \textup{Beta} (\alpha,\beta)+ \mu'-\delta\frac{\alpha}{\alpha+\beta}$ has support
interval of length $\delta$, has mean $\mu'$, and the support interval is contained in the interval $[0,1]$.
\item The distribution $\delta \cdot \textup{Beta} (\alpha,\beta)+ \mu'-\delta\frac{\alpha}{\alpha+\beta}$ is uniform when
$\alpha,\beta$ is $1,1$, and is bimodal when $\alpha,\beta$ is $0.5,0.5$, is symmetric unimodal when $\alpha,\beta$ is $2,2$, and is skewed
in the other two cases, in different directions, one where the mode is an endpoint of the support and one where the mode is interior of the support.
\item The Bernoulli mean $\mu$ is approximately $\mu'$, since ${n+s \choose 2}$ is very large for these
purposes.
\end{itemize}

\newcommand{\figfourtwowidth}{0.45\textwidth}
\begin{figure}[h!]
    \centering
    \begin{subfigure}[b]{\figfourtwowidth}
        \includegraphics[width=1.0\textwidth]{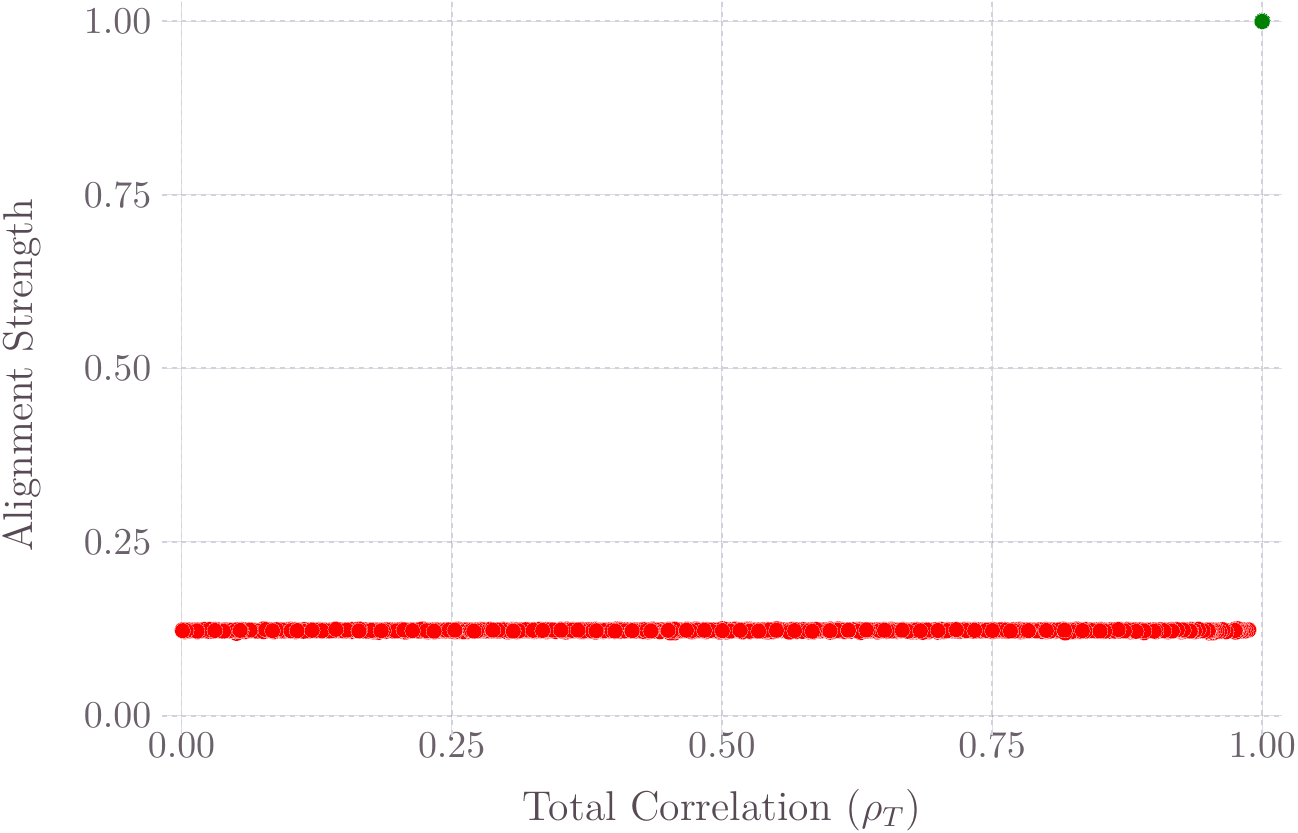}
        \subcaption{number of seeds $s=0$}
        \label{fig:52a:a}
    \end{subfigure}
    \begin{subfigure}[b]{\figfourtwowidth}
        \includegraphics[width=1.0\textwidth]{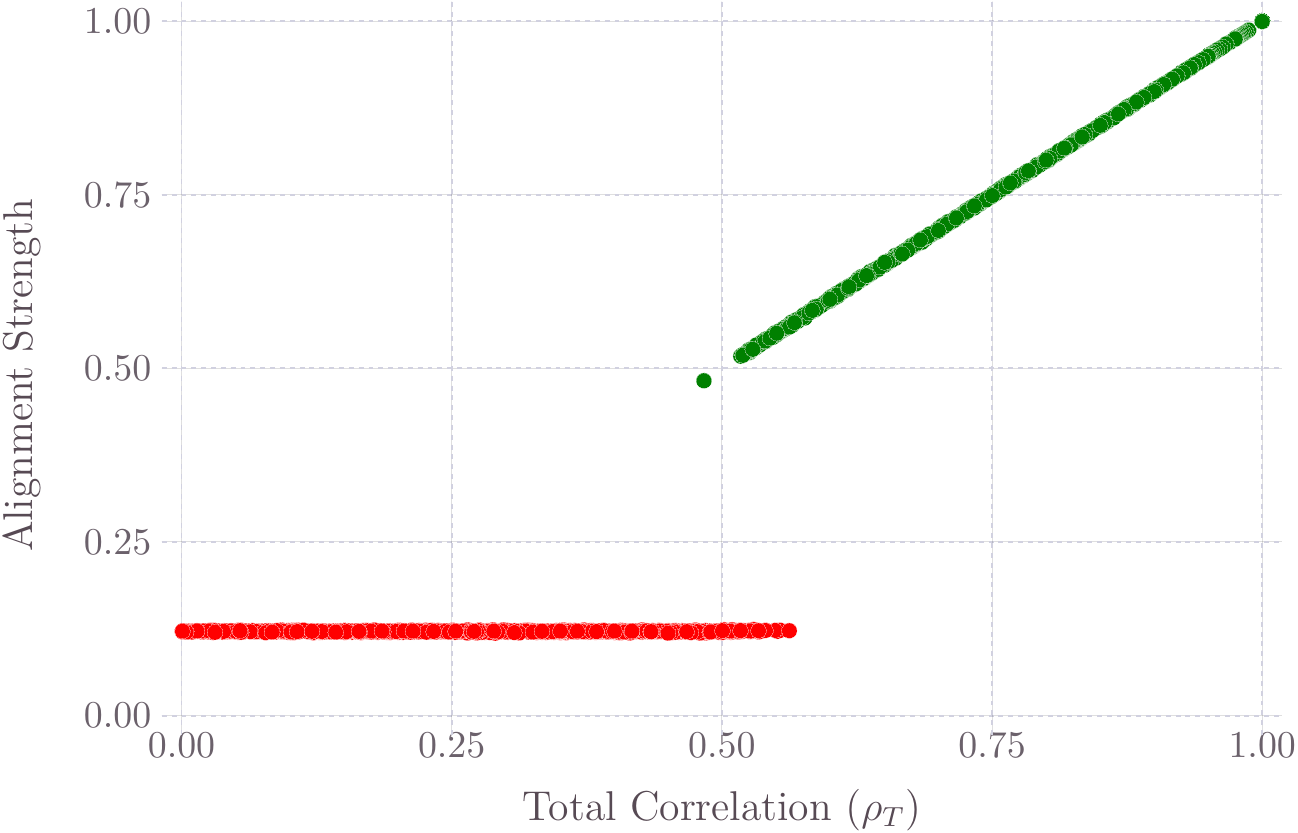}
        \subcaption{number of seeds $s=10$}
        \label{fig:52a:b}
    \end{subfigure}
    \begin{subfigure}[b]{\figfourtwowidth}
        \includegraphics[width=1.0\textwidth]{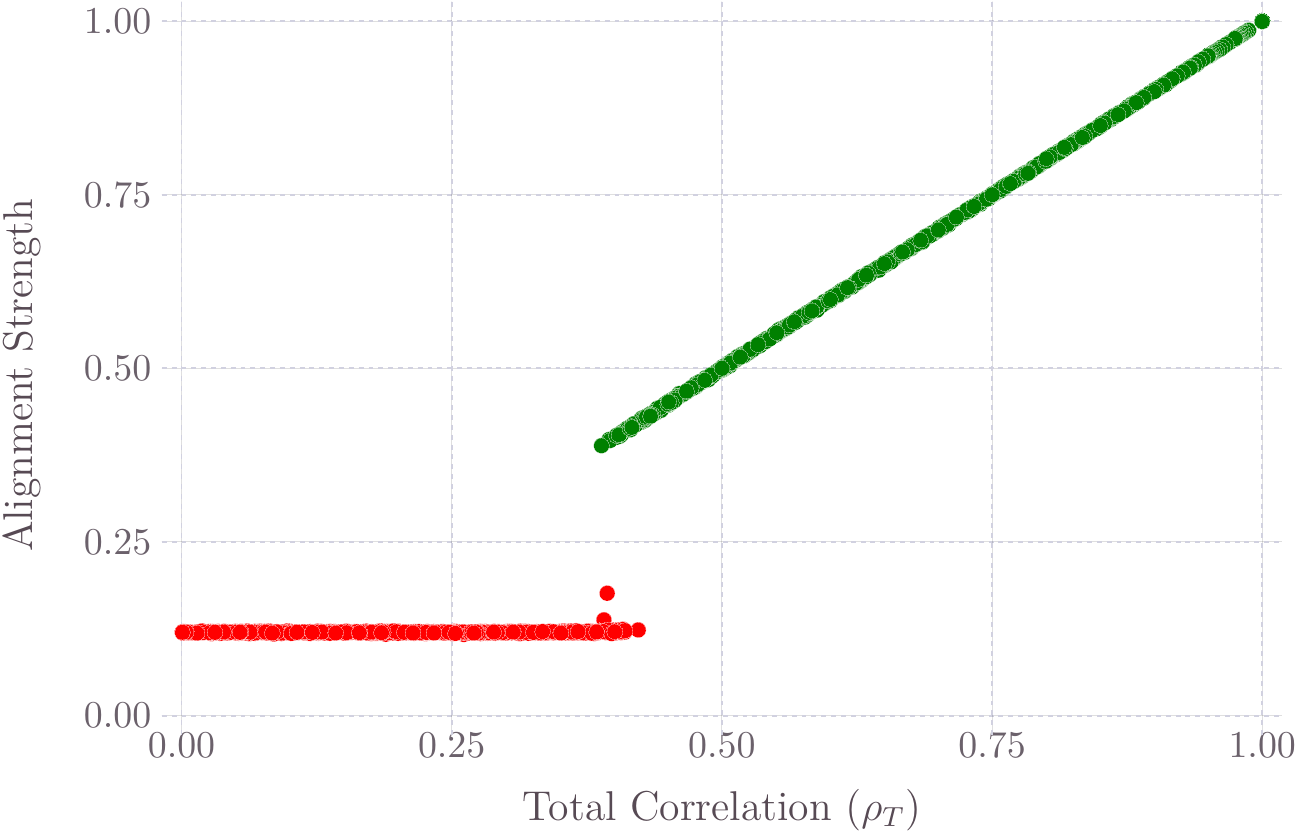}
        \subcaption{number of seeds $s=20$}
        \label{fig:52a:c}
    \end{subfigure}
    \begin{subfigure}[b]{\figfourtwowidth}
        \includegraphics[width=1.0\textwidth]{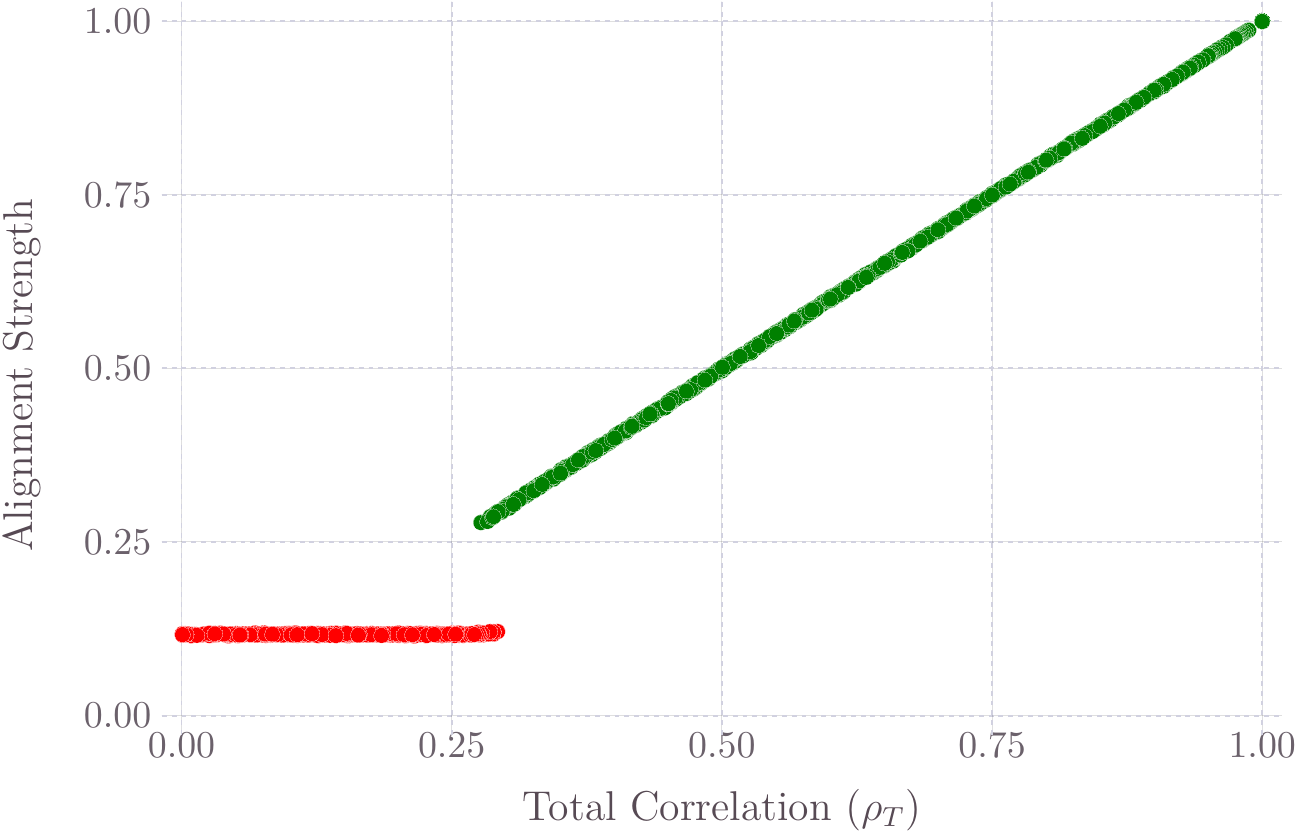}
        \subcaption{number of seeds $s=50$}
        \label{fig:52a:d}
    \end{subfigure}
    \begin{subfigure}[b]{\figfourtwowidth}
        \includegraphics[width=1.0\textwidth]{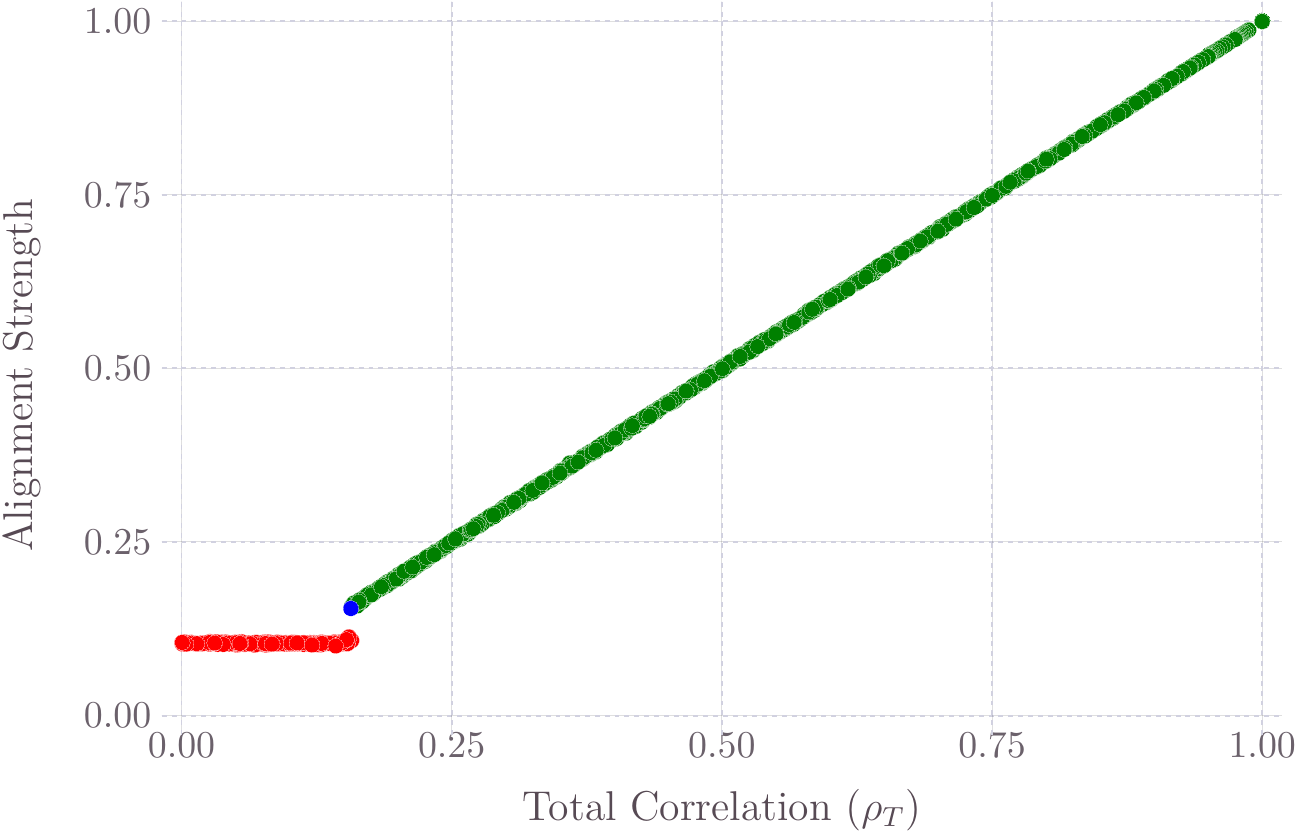}
        \subcaption{number of seeds $s=250$}
        \label{fig:52a:e}
    \end{subfigure}
    \begin{subfigure}[b]{\figfourtwowidth}
        \includegraphics[width=1.0\textwidth]{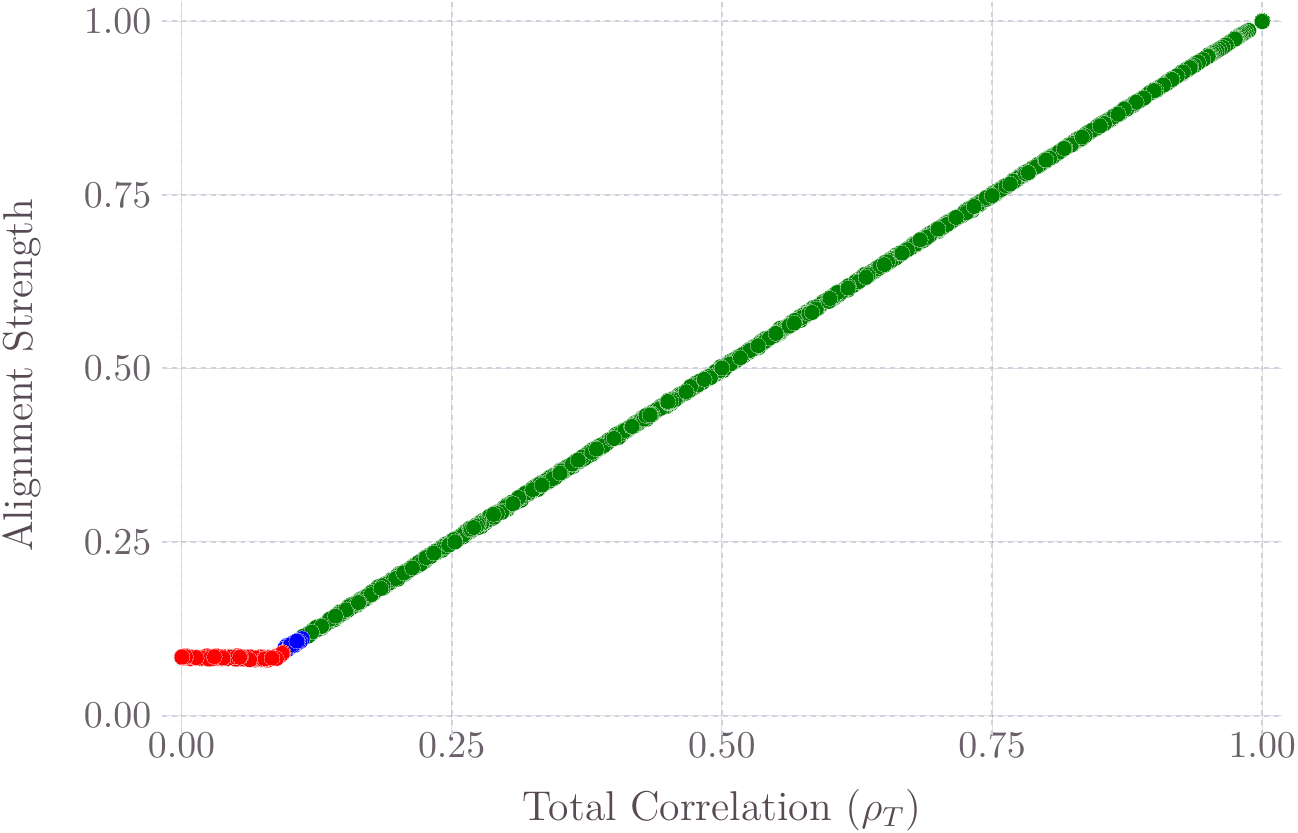}
        \subcaption{number of seeds $s=1000$}
        \label{fig:52a:g}
    \end{subfigure}
      \begin{subfigure}[b]{\figfourtwowidth}
        \includegraphics[width=1.0\textwidth]{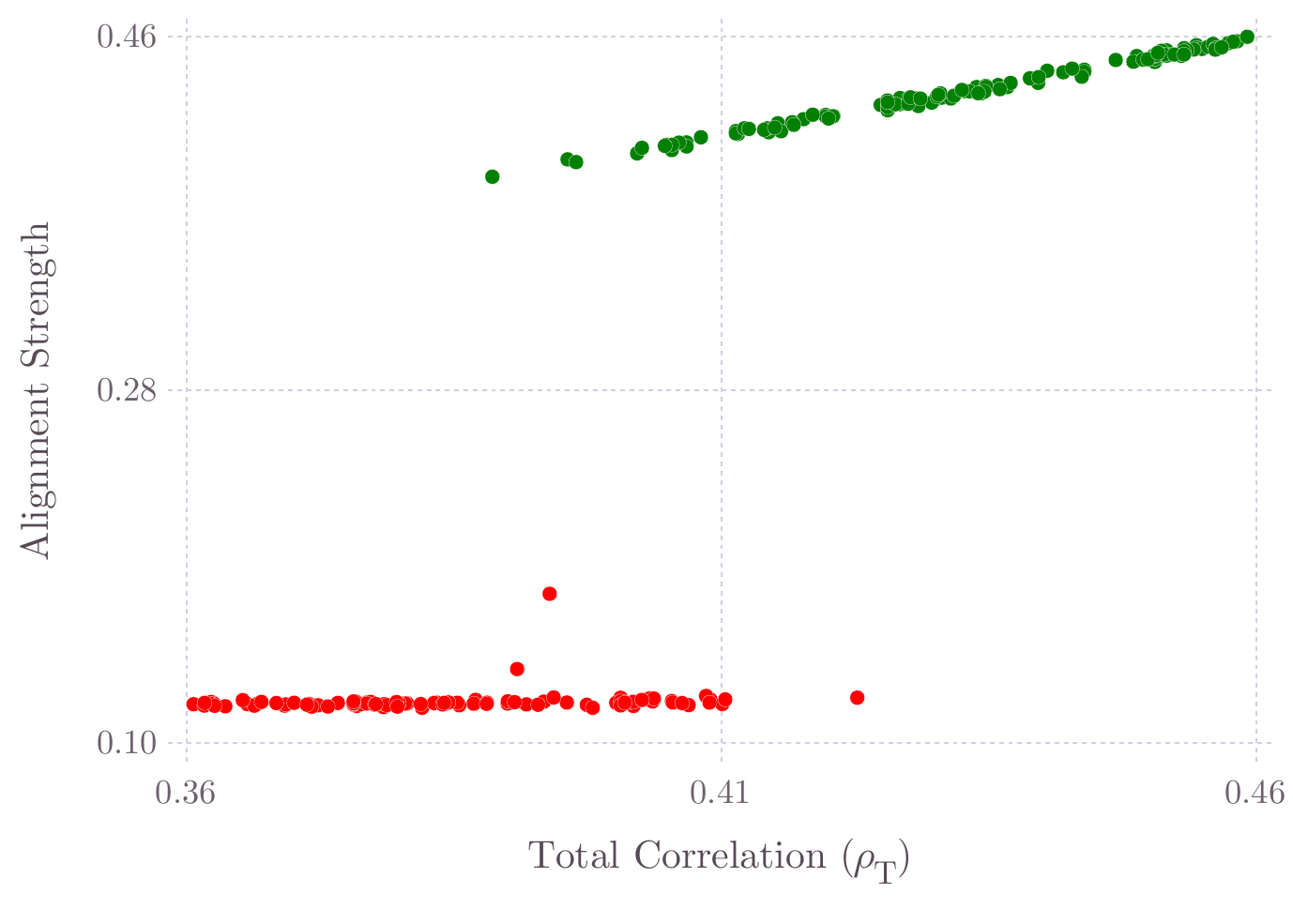}
        \subcaption{Zoom into subfigure (c), $s=20$}
        \label{fig:52a:h}
    \end{subfigure}
      \begin{subfigure}[b]{\figfourtwowidth}
        \includegraphics[width=1.0\textwidth]{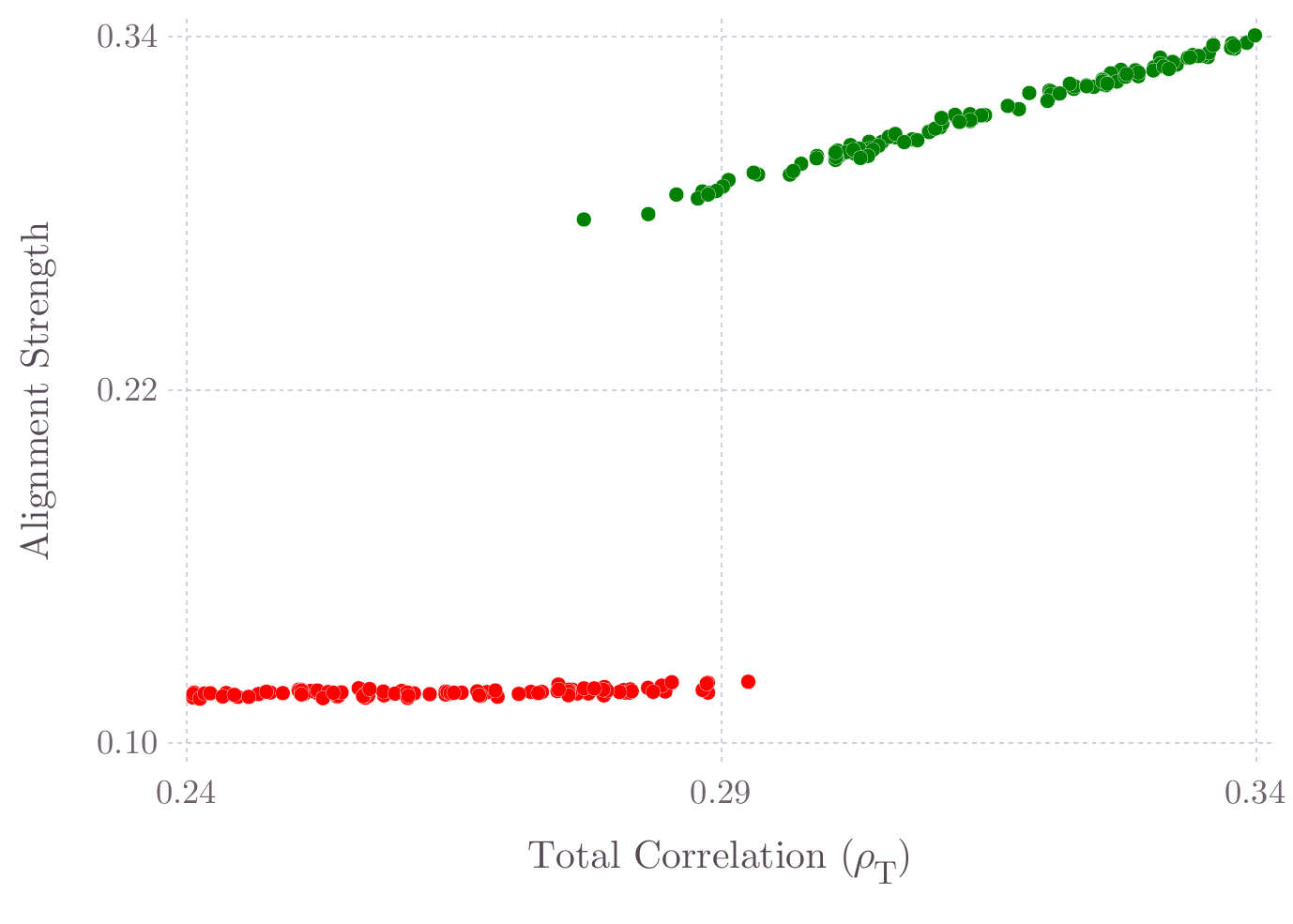}
        \subcaption{Zoom into subfigure (d), $s=50$}
        \label{fig:52a:i}
    \end{subfigure}
      \begin{subfigure}[b]{.2\textwidth}
        \includegraphics[width=1.0\textwidth]{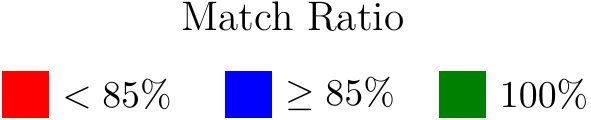}
    \end{subfigure}
    \caption{Alignment strength $\st (\phS)$ plotted against total correlation $\varrho_T$ for the synthetic data experiments in Section \ref{sec:broken}, separated according to the number of seeds $s$. The number of nonseeds  was $n=1000$, and
    only the case of $\mu'=0.5$~is shown here. Match ratio of each experiment is
    color coded green, blue, or red according to the legend above.
    Subfigures (g) and (h) are zooms into subfigures (c) and (d), to increase the granularity so that the thresholding is better seen. }
    \label{fig:52a}
\end{figure}

The $s$ seeds were chosen discrete uniform randomly
from the $n+s$ vertices, and we computed $\phS$ via the SGM algorithm for seeded graph matching.
In Figure \ref{fig:52a}
we plotted alignment strength $\st (\phS)$ against total correlation $\varrho_T$
for all of the pairs of graphs generated in the case where $\mu'=0.5$,
in different subfigures for the different values of $s= 0,10,20,50,250,1000$; green dots indicate
when $\phS=\varphi^*$, blue and red dots indicate when $\phS \ne \varphi^*$, blue when $\phS$ agreed
with $\varphi^*$ on at least $85\%$
of the nonseeded vertices (i.e. ``match ratio $\geq 85\%$''), and red when $\phS$ agreed with $\varphi^*$ on less than $85\%$
of the nonseeded vertices.

Note that in Figure \ref{fig:52a}, each of (a)-(f) are plots of $2255$ points, each point represented with a filled circle,
and the crowding of the points makes them resemble lines; so,
in Figure \ref{fig:52a}, we also included (g) and (h), which are zooms of a portion of (c) and (d), respectively.
With the increased granularity in (g) and (h), we see that if we ignore some outlier red and green dots, then
there is a better defined transition from red to green than would appear in (c) and (d).

The Phantom Alignment Strength Conjecture is well motivated by the results illustrated in
Figure \ref{fig:52a}. In particular, alignment strength $\st (\phS)$ exhibits very low variance
and is approximately a piecewise-linear function of total correlation $\varrho_T$. There appears to be a
critical value  $\qhs$, dependent on the number of seeds $s$ in these experiments, for which the following holds. When
total correlation $\varrho_T$ is above $\qhs$ then $\phS=\varphi^*$
and $\st (\phS) \approx \varrho_T$,
and when total correlation $\varrho_T$ is below $\qhs$ then $\phS \ne \varphi^*$, evidenced by $\st (\phS) \not \approx \varrho_T$, and
$\st (\phS)$ is constant--- at a phantom alignment strength level. When there are enough seeds, we see that  the two pieces of the function
join to become continuous, suggesting that $\phS=\phh$ is then achieved for all $\varrho_T$, and the value of $\qhs$ is then $\qh$.

Also note that the five different Beta distributions
from which Bernoulli parameters were realized (the five pairs of Beta parameters labelled A, B, C, D, E) in these experiments were collected into each of the
figures of Figure \ref{fig:52a}, and
the experiment results for these different distributions
are indistinguishable from each other in the figures,
in accordance with Theorem \ref{thm:new},
and reflected in the Phantom Alignment Strength Conjecture
claim that the phantom alignment strength is just a
function of $n,s,\mu'$, and that it isn't relevant what
distribution is used to obtain the Bernoulli parameters.

Also note the phase transition from matchable to non-matchable which
takes place when $\varrho_T$ gets to $\qhs$, and this phase transition becomes better and better defined as the number
of seeds goes up.

For the other values of $\mu'$, the figures exhibited the same overall type of structure, although the phantom alignment strength values were different. In the interest of space, we
only present here the $\mu'=0.5$ experiment figures.

\subsection{Phantom alignment strength vs theoretical matchability threshold \label{sec:comp}}

Among other assertions,
the Phantom Alignment Strength Conjecture asserts, under
conditions, that the alignment strength $\st (\phS)$
when  $\varrho_T=0$, called the
``phantom alignment strength,'' is equal to the total
correlation threshold for matchability of exact seeded graph matching (i.e. the particular value such that $\phh=\varphi^*$ or not according as $\varrho_T$ is greater than this value or not); indeed, we have denoted this common quantity~$\qh$.
In this section, we will compare alignment strength $\st (\phS)$ when  $\varrho_T=0$
to the matchability threshold proved in
\cite{JMLR:v15:lyzinski14a}.

Consider a probability space with a sequence of correlated Bernoulli random graphs for each of the number of vertices $\n \equiv n=1,2,3,\ldots$, with $s=0$ seeds and all Bernoulli parameters
equal to a fixed value $p$ (ie correlated Erdos-Renyi random graphs).
When we say that a sequence of events happens ``almost always'' we mean that, with probability $1$, all but a finite number of
the events occur.
The following result was stated and proved in \cite{JMLR:v15:lyzinski14a}; although stated there in terms of $\varrho_e$, we write $\varrho_T$ instead,
since here, where $\varrho_h=0$, we have that $\varrho_T=\varrho_e$.

\begin{theorem} There exists positive, real valued, fixed constants $c_1,c_2$ such that if
$\varrho_T \geq c_1 \sqrt{\frac{\log n}{n}}$ then almost always $\phh=\varphi^*$, and
if $\varrho_T \leq c_2 \sqrt{\frac{\log n}{n}}$ then
$\lim_{n \rightarrow \infty} \e | \{  \varphi \in \varPi : \dis' (\varphi)< \dis'( \varphi^*  ) \}  |=\infty$.
\end{theorem}

\begin{figure}[h!]
    \centering
    \includegraphics[width=0.7\textwidth]{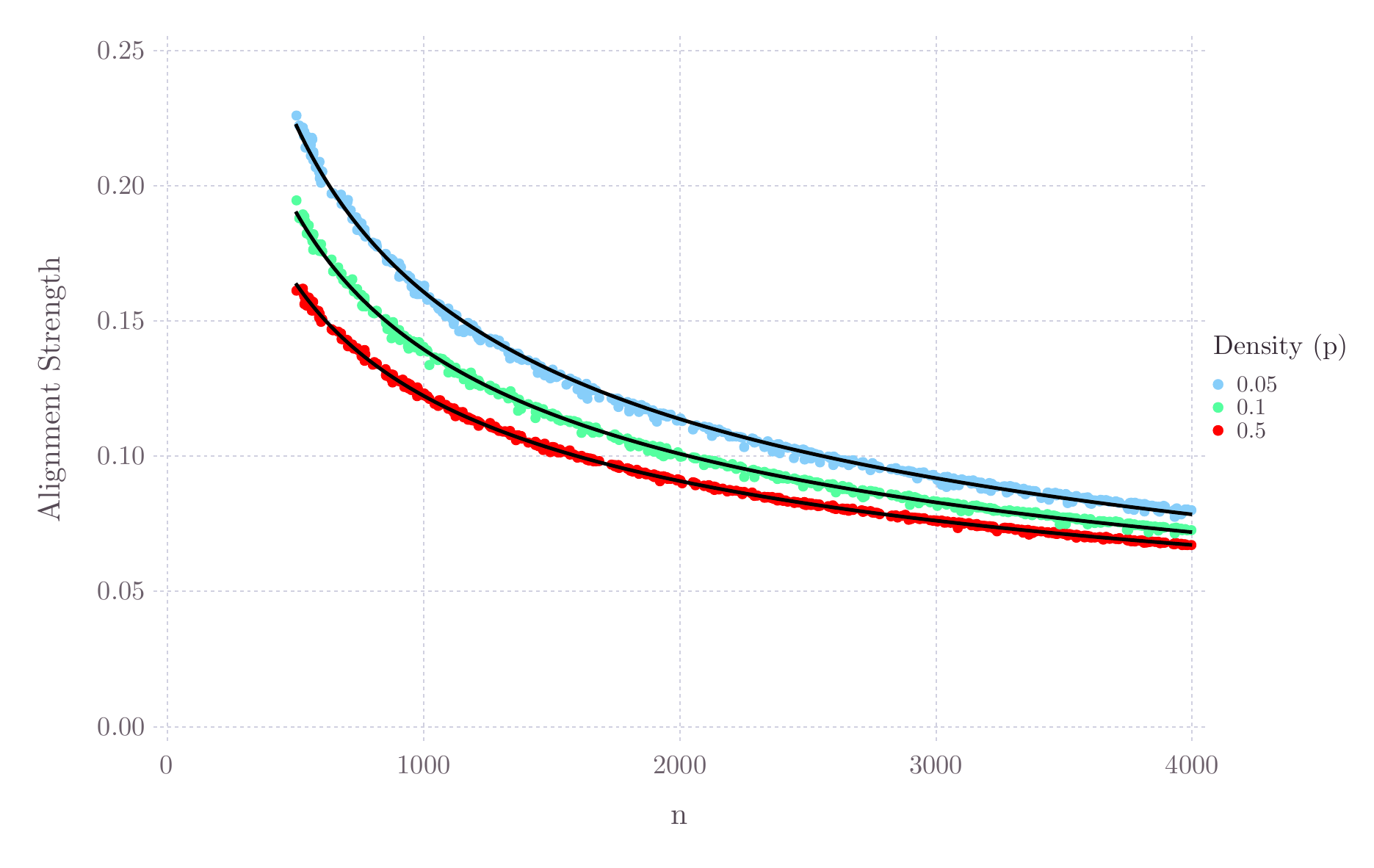}
    \caption{Phantom alignment strength as a function of $n$,\\ fitted to $f_p(n):=d_p+ c_p \sqrt{\frac{\log n}{n}}$.}
    \label{fig:53a}
\end{figure}

For each value of $p=.05, .1, .2, .3, .4, .5$, and
each of $500$ values of $n$ between $500$ and $4000$, (as mentioned, $s=0$) we plotted realizations of alignment strength $\st (\phS)$ vs the value of~$n$, for uncorrelated ($\varrho_e=0$) pairs of random Bernoulli (Erdos-Renyi) graphs where each Bernoulli parameter is~$p$, hence $\varrho_T=0$ (since $\varrho_e=0$, $\varrho_h=0$). Figure \ref{fig:53a} shows the plots for $p=0.05,0.1,0.5$.

 Then, for each $p$, we fit the associated points to a curve $f_p(n):=d_p+ c_p \sqrt{\frac{\log n}{n}}$ for real numbers $c_p$ and $d_p$;
 the values of $d_p$ and $c_p$ are given in Table \ref{table:53a}, and $f_p$ is also drawn in Figure \ref{fig:53a}. For each
value of $p$, note the near-perfect fit of $f_p$ to the
associated points plotted in Figure \ref{fig:53a}, and note that the value of $d_p$ is close to zero.

\begin{table}[h!]
    \centering
    \caption{Values of the constants in $f_p(n):=d_p+ c_p \sqrt{\frac{\log n}{n}}$ }
    \label{table:53a}
    \begin{tabular}{l|rr}
        \hline\hline
        $p$ & $d_p$ & $c_p$ \\
        \hline
        $0.05$ &    $-0.021$    & $2.19$ \\
        $0.1$  &	$-0.010$    & $1.80$ \\
        $0.2$  &	$-0.003$    & $1.58$ \\
        $0.3$  &	$-0.001$    & $1.51$ \\
        $0.4$  &	$ 0.000$    & $1.48$ \\
        $0.5$  &	$ 0.000$    & $1.47$
    \end{tabular}
\end{table}

Indeed, this suggests, as conjectured in the Phantom
Alignment Strength Conjecture,
that the phantom alignment strength (ie $\st (\phS)$
when  $\varrho_T=0$) exists as a value $\qh$ which coincides with  the amount of total correlation needed for $\phh=\varphi^*$.

\subsection{Block settings \label{sec:block}}

The setting of the Phantom Alignment Strength Conjecture in Section \ref{sec:ASC} was specifically concerning
correlated Bernoulli random graphs $G_1,G_2$ such that there are $n$ nonseed vertices,
$s$ seed vertices (selected discrete uniformly from the $\n:=n+s$ vertices), and Bernoulli parameters for each
pair of vertices are selected independently from any fixed distribution with mean $\mu'$.

Let us consider a block setting, which differs from the above in that there is a positive integer~$K$, and the
vertex set $V$ is first randomly partitioned into $K$ blocks $B_1,B_2,\ldots,B_K$ as follows:
There is a given probability vector $\pi \in [0,1]^K$ such that $\sum_{i=1}^K\pi_i=1$ and
each vertex in $V$ is independently placed in block $B_i$ with probability $\pi_i$ for $i=1,2,\ldots,K$.
Next, suppose there is a unit-interval-valued (ie $[0,1]$-valued)
distribution $F_{i,j}$ for each $i=1,2,\ldots,K$ and $j=i,i+1,\ldots,K$ such that, for each
$1\leq i \leq j \leq K$ and each $u\in B_i$ and $v \in B_j$, the Bernoulli parameter $p_{\{ u,v \}}$
is independently realized from distribution $F_{i,j}$. Let $M$ be the $K\times K$ symmetric matrix with $i,j$th
entry equal to the mean of distribution $F_{i,j}$.

Similarly to the Phantom Alignment Strength Conjecture, does there exists a phantom alignment strength value
$\qh \equiv \qh (n,s,\pi,M) \in [0,1]$ and also $\qhs \equiv \qhs (n,s,\pi,M) \in [0,1]$ whereby
Equation (\ref{eqn:one}) and Equation (\ref{eqn:two}) hold? This is not so simple.

\begin{figure}[h!]
    \centering
    \includegraphics[width=0.65\textwidth]{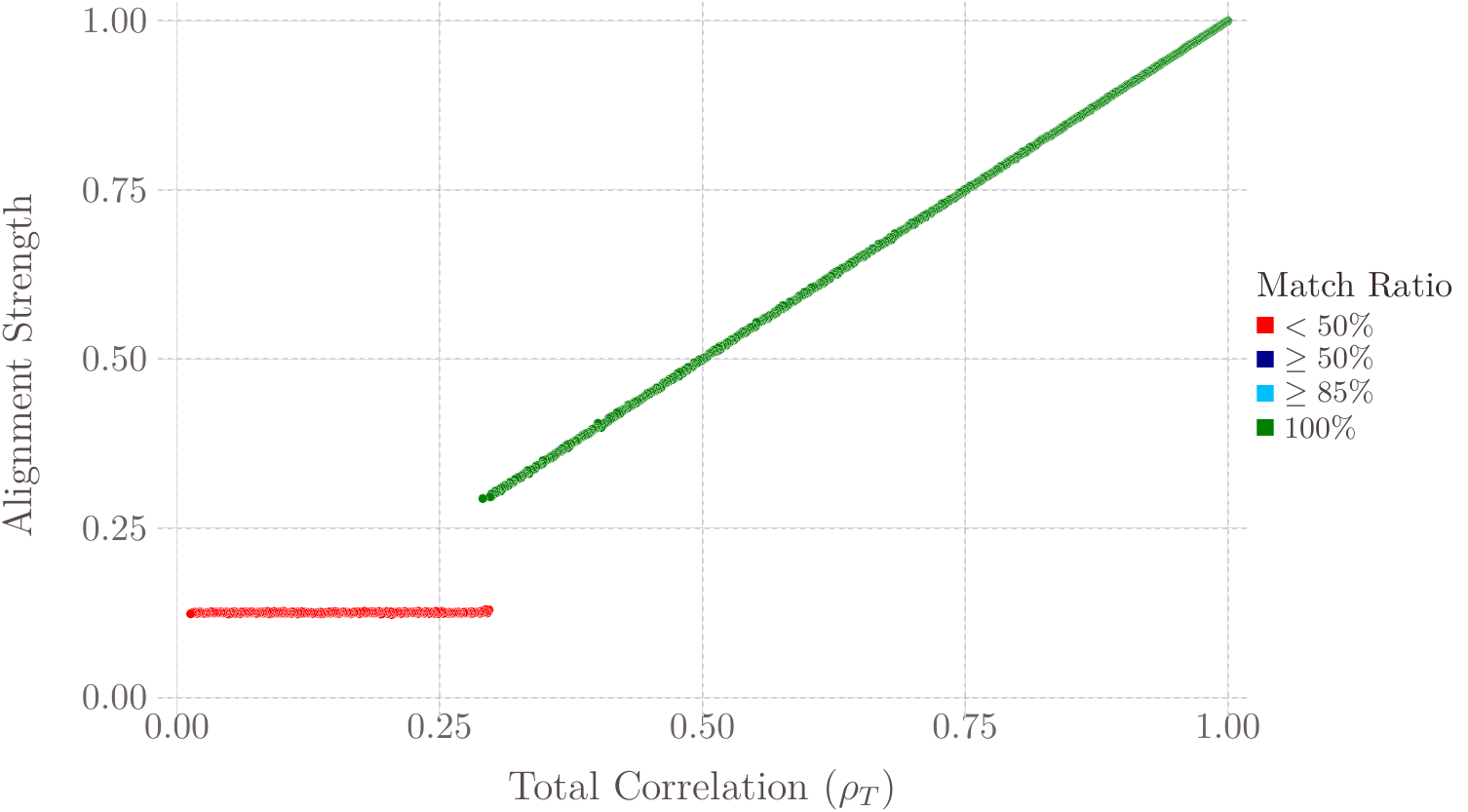}
    \caption{Experiment A in Section \ref{sec:block}; here $F_{1,1}$, $F_{1,2}$, $F_{2,2}$ are resp. point mass at $0.3$, $0.4$, $0.5$. }
    \label{fig:54a}
    \centering
    \includegraphics[width=0.65\textwidth]{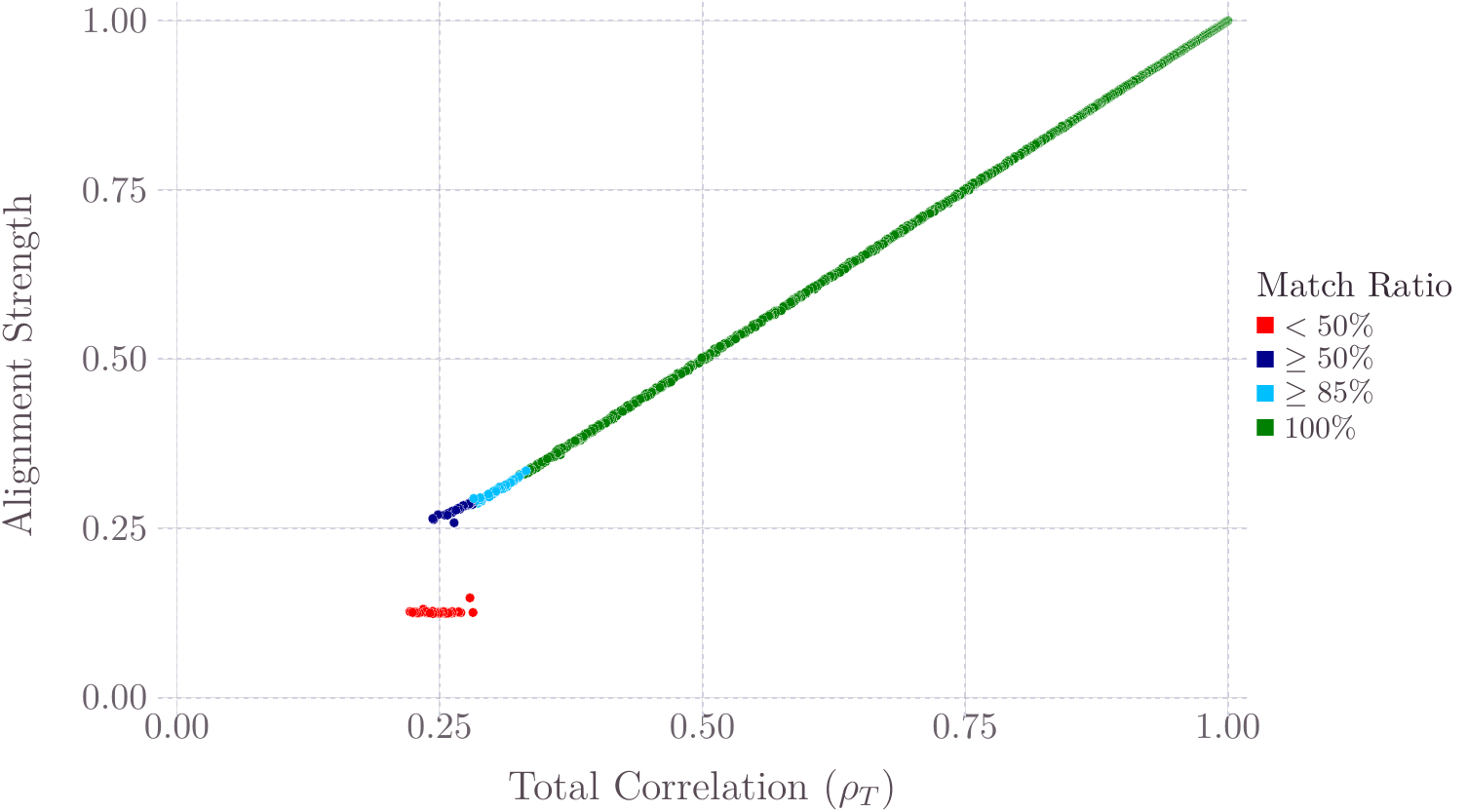}
    \caption{Experiment B in Section \ref{sec:block}; same as
    Experiment A except that $F_{2,2}$ is  uniform  $[0,1]$.}
    \label{fig:54b}
    \centering
    \includegraphics[width=0.65\textwidth]{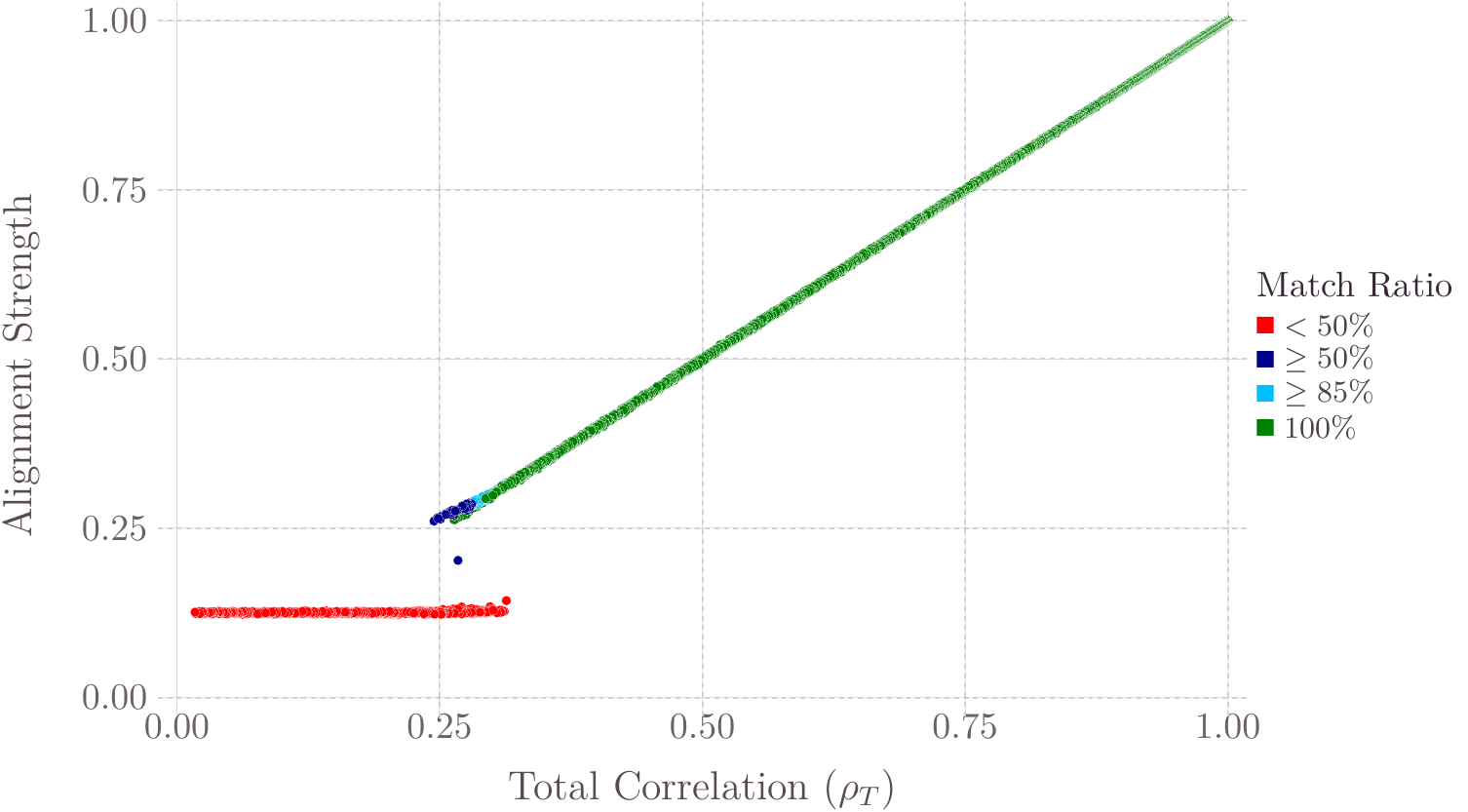}
    \caption{Experiment C in Section \ref{sec:block}; eight different combinations for $F_{1,1}$, $F_{1,2}$, $F_{2,2}$. }
    \label{fig:54c}
\end{figure}

We consider the following choices for $n$, $s$, $\pi$, and $M$:
\begin{eqnarray*}
n=1000 \ \ \ \  s=40 \ \ \ \
\pi = \left [  \begin{array}{c} 0.2 \\ 0.8 \end{array} \right ] \ \ \ \
M= \left [  \begin{array}{cc}  0.3 & 0.4 \\  0.4 & 0.5 \end{array} \right ]
\end{eqnarray*}
In experiment ``A'', we took $F_{1,1}$ to be point mass distribution at $0.3$, $F_{1,2}$ to be point mass distribution at $0.4$, and
$F_{2,2}$ to be point mass distribution at $0.5$. For each value of edge correlation $\varrho_e$ from $0$ to $1$
in increments of $0.001$, we realized Bernoulli parameters
and then we realized associated correlated Bernoulli random graphs.
In Figure \ref{fig:54a}, we plotted alignment strength $\st (\phS)$ against total correlation $\varrho_T$;
green dots indicate when $\phS=\varphi^*$, (else)
light blue when $\phS$ agreed with $\varphi^*$ on at least $85\%$ of the nonseeded vertices, (else)
dark blue when $\phS$ agreed with $\varphi^*$ on at least $50\%$, (else) red
when $\phS$ agreed with $\varphi^*$ on less than $50\%$ of the nonseeded vertices. We then repeated the
experiment with the only difference being that $F_{2,2}$ was the uniform distribution on the interval $[0,1]$,
so $(n,s,\pi,M)$ are same as above; the
resulting plot is Figure~\ref{fig:54b} (alignment strength $\st (\phS)$ vs $\varrho_T$, same dot color scheme as above). Let us call this Experiment ``B.''

Next, we repeated the above experiment for all eight possible combinations of:\\
$F_{1,1}$ is the uniform distribution on a) interval $[0.25,0.35]$ or b) interval $[0,0.6]$\\
$F_{1,2}$ is the uniform distribution on a) interval $[0.35,0.45]$ or b) interval $[0,0.8]$\\
$F_{2,2}$ is the uniform distribution on a) interval $[0.45,0.55]$ or b) interval $[0,1]$\\
and we superimposed all of the alignment strength vs total correlation plots in Figure \ref{fig:54c} (same dot color scheme as above); we will call this Experiment ``C.''  Again, the underlying $(n,s,\pi,M)$ are the same as the previous experiments.

Note that Figure \ref{fig:54a}, Figure  \ref{fig:54b}, and
Figure \ref{fig:54c} (for respective experiments A,B, and C) are not similar, even though they originate
from the same values of $n$, $s$, $\pi$, and $M$. Thus, the
Phantom Alignment Strength Conjecture is not simply
extended to the case of nontrivial block structure.

However, also note that when SGM was broadly failing to get the truth in experiments A, B, and C (i.e. the red dots in
Figure \ref{fig:54a}, Figure  \ref{fig:54b}, and
Figure \ref{fig:54c}), the alignment strength was almost
constant, at a value of around $0.12$. This suggests a decision  procedure (analogous the procedure
described in Section \ref{sec:ASC}) for deciding if $G_1$, $G_2$ from an
$(n,s,\pi,M)$-block model are graph matched with some truth.
The procedure would be to realize $H_1$ and $H_2$ as correlated Bernoulli random graphs where $\varrho_e=0$,
where the $n+s$ vertices are apportioned to the blocks
in proportion to $\pi$, and where, for every pair of vertices,
the Bernoulli parameter is taken as the entry of $M$
associated with the block memberships of the two vertices,
and then the $s$ seeds are chosen uniformly at random.
The alignment strength of
the seeded graph match of $H_1$ to $H_2$ can then be used as
a phantom alignment strength value in the sense that, if the alignment strength
of the seeded graph match of $G_1$ to $G_2$ is more than some $\epsilon>0$ greater than this phantom
alignment strength value, then we decide that there is
at least some truth present in the seeded graph match of $G_1$
to $G_2$.

What made the block structure more complicated? We will next provide some insight. Indeed,
Experiment~B was constructed in an
extreme way in order to cause particular mischief.
The value of $\varrho_h$ in Experiment~A was
approximately $.0129$, and the value of
$\varrho_h$ in Experiment B was approximately $.2277$; in
particular, that is why the value of $\varrho_T$ was never
below approximately~$.22$ in Experiment B, as is clear from Figure  \ref{fig:54b}. However, in Experiment B when
$\rho_e=0$, all of the vertices in the first block are stochastic twins; they share the same probabilities
of adjacency as each other to all of the vertices in the graph, and all adjacencies are collectively independent.
Thus the ``true'' bijection (the identity) has no
signal in that case. (One might even say that the ``truth''
isn't very ``truthy.'')
As such, the total correlation in that case, approximately
$.2277$, does not contribute to matchability vis-a-vis the
first block. As positive edge correlation $\varrho_e$ is
increasingly added in to Experiment~B, the first block
achieves matchability on the strength of only the edge
correlation, and the second block achieves matchability
on the strength of edge correlation together with
heterogeneity correlation. In this manner, total
correlation does not tell a uniform story across all
vertices. This is in contrast to the hypotheses of
the Phantom Alignment Strength Conjecture (and the
setup in the empirical matchability experiments in the paper \cite{fishkind2019alignment})
where the Bernoulli parameters were realized from one distribution. Note that with Experiment C, there is more
variety in $\varrho_h$ (for the eight experiments the values of $\varrho_h$ ranged from approximately $.0161$ to approximately $.30$); there is still some
lack of demarcation
between matchable and nonmatchable in terms of total correlation, but the situation is improved somewhat from
the left tail of the figure, and total correlation has more influence as a unified quantity.

We did additional experiments with other values
of $(n,s,\pi,M)$ and found comparable results to what
appears above.

\subsection{Real data; matching graphs to noisy renditions  \label{sec:noisy}}

Recall that the Phantom Alignment Strength Conjecture
is formulated under the assumption that each
pair of vertices has a Bernoulli parameter that
is a realization of a distribution which is common
to all of the pairs of vertices. How
realistic is this assumption in practice? And, more
to the point of the practitioner, do the conclusions
of the conjecture apply to real data, in general?

In this section we consider a human connectome at
different resolution levels. (This connectome has been
featured in \cite{Priebe2019Mar,Chung2020Aug}.)
Diffusion-weighted Magnetic Resonance Imaging (dMRI) brain scans were collected from one hundred and fourteen humans at the Beijing Normal University~\cite{corr}. Fiber tracts, which trace axonal pathways through a three-spatial-dimensional cuboid array of $1\times 1 \times 1 \;\textrm{mm}^3$  \textit{voxels} of the dMRI scan, are estimated using the ndmg pipeline \cite{ndmg}.

For each value of $\n= 70, \ 107, \ 277, \ 582, \ 3230$, the graph $G_{\n}$
was formed in the following manner.
Starting from the original cuboid array
of voxels, $\n$ equally spaced ``contractile'' voxels were selected, and each voxel in the array
was merged with its nearest contractile voxel \cite{glocal};
the $\n$ such groupings of voxels (centered at their
contractile voxel) are the $\n$ vertices of
the graph~$G_{\n}$. For any two vertices in $G_{\n}$, we declare them
adjacent precisely when there exists a fiber that runs
through any voxel of one vertex and also any voxel of
the other vertex for any of the one hundred and fourteen individuals.

Given any graph $G=(V,E)$,
and also given any noise parameter
$\rho \in [0,1]$, we can instantiate a
graph $\overline{G}$ called a
{\it $\rho$-noised rendition of} $G$ on the same
vertex set $V$ as follows.
Denote the density of $G$
by $\dens'G:=\frac{|E|}{{|V| \choose 2}}$. First,
instantiate an independent Erdos-Renyi graph $H$ on $V$ with
Bernoulli parameter $\dens'G$; i.e.~each
pair of vertices is an edge independently of the others with probability~$\dens' G$.
Next, for each pair of vertices $\{u,v\}$,
perform an independent Bernoulli trial; with probability
$\rho$ set $u$ adjacent/ not adjacent (resp.) to $v$ in $\overline{G}$ according as $u$ adjacent/ not adjacent (resp.) to $v$ in $G$, and
with probability $1-\rho$ set
$u$ adjacent/ not adjacent (resp.) to $v$ in $\overline{G}$ according as $u$ adjacent/ not adjacent (resp.) to $v$ in $H$.
In this manner, $\overline{G}$ is a mixture of $G$ and noise graph $H$.
When graph matching $G$ to a $\rho$-noised rendition of $G$,
clearly $\varphi^*$ is the identity function $V$ to $V$.

\begin{figure}[h!]
    \centering
    \includegraphics[width=1.0\textwidth]{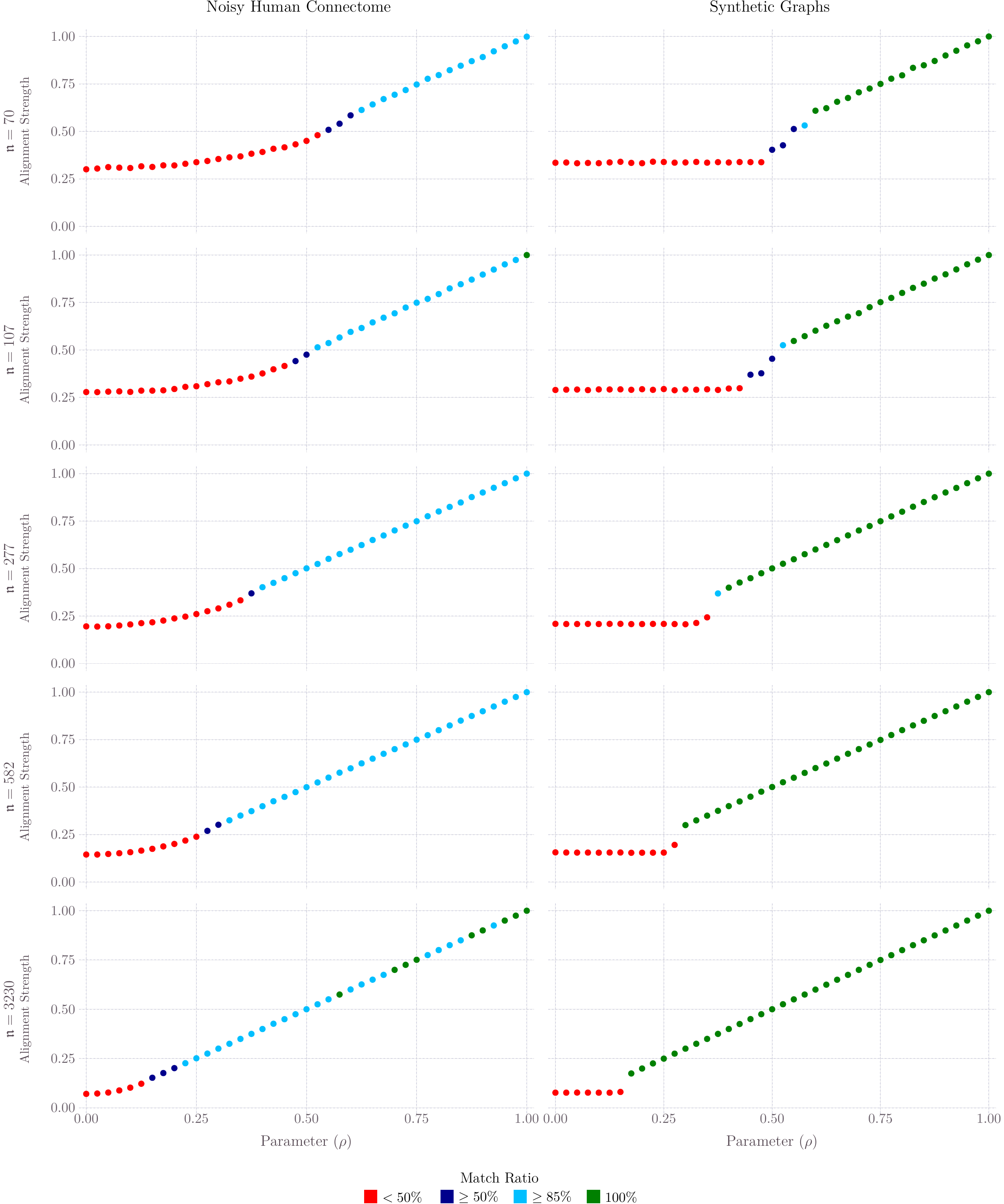}
    \caption{Section \ref{sec:noisy} experiments; LHS is noisy connectome, RHS is corresponding synthetic.}
    \label{fig:55}
\end{figure}

For each of $\n= 70, \ 107, \ 277, \ 582, \ 3230$, we did the following experiment. For each value of the noise parameter $\rho$ from $0$ to $1$ in increments of $.025$, we did $20$ repetitions of
instantiating a $\rho$-noised rendition of $G_{\n}$,
then seeded graph matched $G_{\n}$
to it using the SGM algorithm after selecting $10\%$
of the $\n$ vertices (discrete uniform randomly) as seeds.
The mean
alignment strength $\st (\phS)$ (the mean being over the $20$
repetitions) vs noise parameter $\rho$ was plotted in five respective figures (for the five different values of $\n$)
in the left side of Figure \ref{fig:55};
green dots indicate when $\phS=\varphi^*$, (else)
light blue when $\phS$ agreed with $\varphi^*$ on at least $85\%$ of the nonseeded vertices, (else)
dark blue when $\phS$ agreed with $\varphi^*$ on at least $50\%$, (else) red
when $\phS$ agreed with $\varphi^*$ on less than $50\%$ of the nonseeded vertices.

We then repeated the above experiments, with the only difference being that in place of $G_{\n}$ we used an
Erdos-Renyi graph instantiation, the Erdos-Renyi using  the Bernoulli parameter $\dens'G_{\n}$ (the density of the connectome $G_{\n}$). The resulting plots are in the
right hand side of Figure \ref{fig:55}. Simple calculations
of the distributions
show that the pairs of graphs being seeded graph matched here in these repeated experiments
are precisely correlated Erdos-Renyi graphs with the
parameter $\rho$ being precisely the edge correlation
$\varrho_e$, which is equal to $\varrho_T$ since $\varrho_h=0$.

To emphasize: The left hand side of Figure  \ref{fig:55} is
from seeded graph  matching connectome to noisy connectome, and the right
hand side of Figure  \ref{fig:55} is from seeded graph matching
synthetic data of the same connectome density to a noisy
version of this synthetic data, which turns out to precisely
be seeded graph matching pairs of correlated Bernoulli random graphs where the noise parameter turns out to be the total
correlation,
so the figures in the right hand side of Figure  \ref{fig:55} are of Section \ref{sec:broken} variety (except that
the alignment strength values are averaged over $20$ instantiations).

Notice that the figures in the left hand side of Figure
\ref{fig:55} and their respective counterparts in the
right hand side  of Figure
\ref{fig:55} look remarkably similar in many important
ways. The differences seem to just be that the
seeded graph matching success and alignment strength values
clearly exhibit thresholding in the synthetic data,
which is less pronounced and more gradual
in the connectome data, although the sharpness
of the connectome thresholding seems to
be catching up as the number of vertices increases.
Aside from this, there stills seems to be a reasonable phantom
alignnment strength for the connectome data.

\subsection{Real data; matching same objects under different modalities \label{sec:modalities}}

In this section, we illustrate the ideas in this paper using
three real data sets from \cite{FAP}; they are
the Wikipedia, Enron, and C Elegans pairs of graphs.
Each is an example of a pair of graphs with the same
underlying objects (thus there is a natural
``true'' bijection), and adjacencies between objects in the
respective graphs are relationships among
the objects in two different modalities.

The Wikipedia pair of graphs $G_1$, $G_2$ from \cite{FAP}
were created in the year 2009.
The vertices of $G_1$ are the English language
Wikipedia articles hyperlinked from
the Wikipedia article ``Algebraic Geometry,'' and all
Wikipedia articles hyperlinked from these articles;
in total, there are $\n=1382$ vertices.
These vertices/articles each
have directly
corresponding articles in the French language Wikipedia,
and these are the vertices of $G_2$.
Every pair of vertices/articles in $G_1$
are adjacent in $G_1$ precisely when one of the articles
hyperlinks to the other article in the English language
Wikipedia, and every pair of
vertices in $G_2$ are adjacent in $G_2$ precisely when
one of the articles links to the other
in the French language Wikipedia.
Thus $G_1$ and $G_2$ are
simple, undirected graphs, and the ``true'' bijection
is the function mapping English articles to their
French versions.

\begin{figure}[h!]
    \centering
    \includegraphics[width=0.6\textwidth]{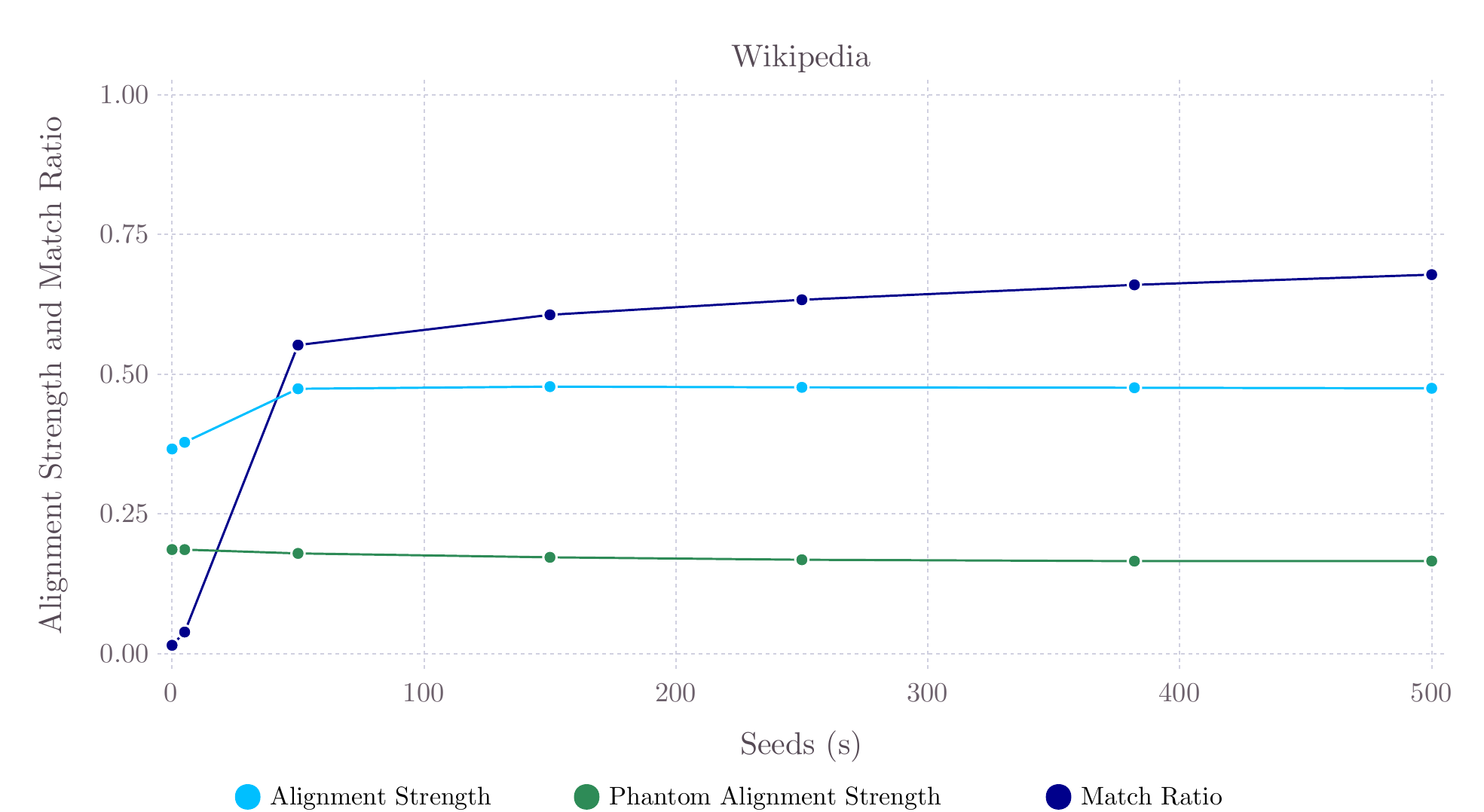}
    \caption{Matching English and French Wikipedia graphs}
    \label{fig:56:wikipedia}
\end{figure}

For each value of $s=0, 5, 50, 150, 250, 382, 500$, we did 100 replicates
of uniformly sampling $s$ seeds from the $\n$ vertices,
seeded graph matched $G_1$ to $G_2$ using SGM, then
recording the
alignment strength $\st (\phS)$, averaged over the 100 replicates, plotted (in blue)
vs the number of seeds $s$ in Figure \ref{fig:56:wikipedia}.
In the same figure, we recorded the match ratio (the
number of nonseeds correctly matched, divided by the
number of nonseeds), averaged over the 100 replications,
plotted (in purple) vs the number of seeds $s$,
also in Figure \ref{fig:56:wikipedia}.
In addition, for each value of $s=0, 5, 50, 150, 250, 382, 500$,
we did 100 replicates of realizing uncorrelated pairs
of Erdos-Renyi graphs $H_1$, $H_2$, each
with $1382$ vertices and
Bernoulli parameter of $H_1$ equal to the density of $G_1$,
Bernoulli parameter of $H_2$ equal to the density of $G_2$,
then uniformly sampling $s$ seeds from the $1382$ vertices,
then seeded graph matched $H_1$ to $H_2$ using SGM, and
recording the alignment strength $\st (\phS)$,
averaged over the 100 replicates, plotted (in green)
vs the number of seeds $s$ in Figure \ref{fig:56:wikipedia}; these values
represent the phantom alignment strength values in the
respective seed levels. Note that, as the number
seeds went from $0$ to $5$ to $50$, the jump
in match ratio coincides with a jump in the gap
between seeded graph matching alignment strength and the phantom alignment strength. (Even when $s=0$ there is some truth in the graph match; the match ratio was $.0151$,
approximately $21$ nonseed vertices matched correctly,
whereas chance is $1/1382$, one nonseed vertex matched correctly.)

The C. Elegans pair of graphs $G^{e \ell}$, $G^{ch}$
from  \cite{Varshney2011Feb,FAP} are connectomes mapping out the
neural structure of the roundworm {\it Caenorhabditis Elegans}.
C. Elegans is of interest to neuroscientists due to
its well studied genetics \cite{Celegans1998genome}, comparatively
simple nervous system \cite{White1986Nov}, and a growing understanding
of the correspondence between the two \cite{Bargmann1998Dec,xn2018Feb}.
Like in humans, communication in the C. Elegans nervous system occurs
via synapses, or junctions, between pairs of neurons.
Neuronal synapses in the C. Elegans connectome can be classified in
two ways \cite{Varshney2011Feb}:
an electrical synapse is a channel through which electrical
impulses traverse,
whereas chemical synapses are junctions through which
neurotransmitters flow.
We consider $\n = 279$ somatic neurons
of the hermaphrodite C. Elegans
as the vertices
of each graph.
For each pair of vertices/neurons,
they are adjacent in $G^{e \ell}$ precisely when there
is an electrical synapse between them, and they are
adjacent in $G^{ch}$ precisely when there is a chemical
synapse between them.

\begin{figure}[h!]
    \centering
    \includegraphics[width=0.6\textwidth]{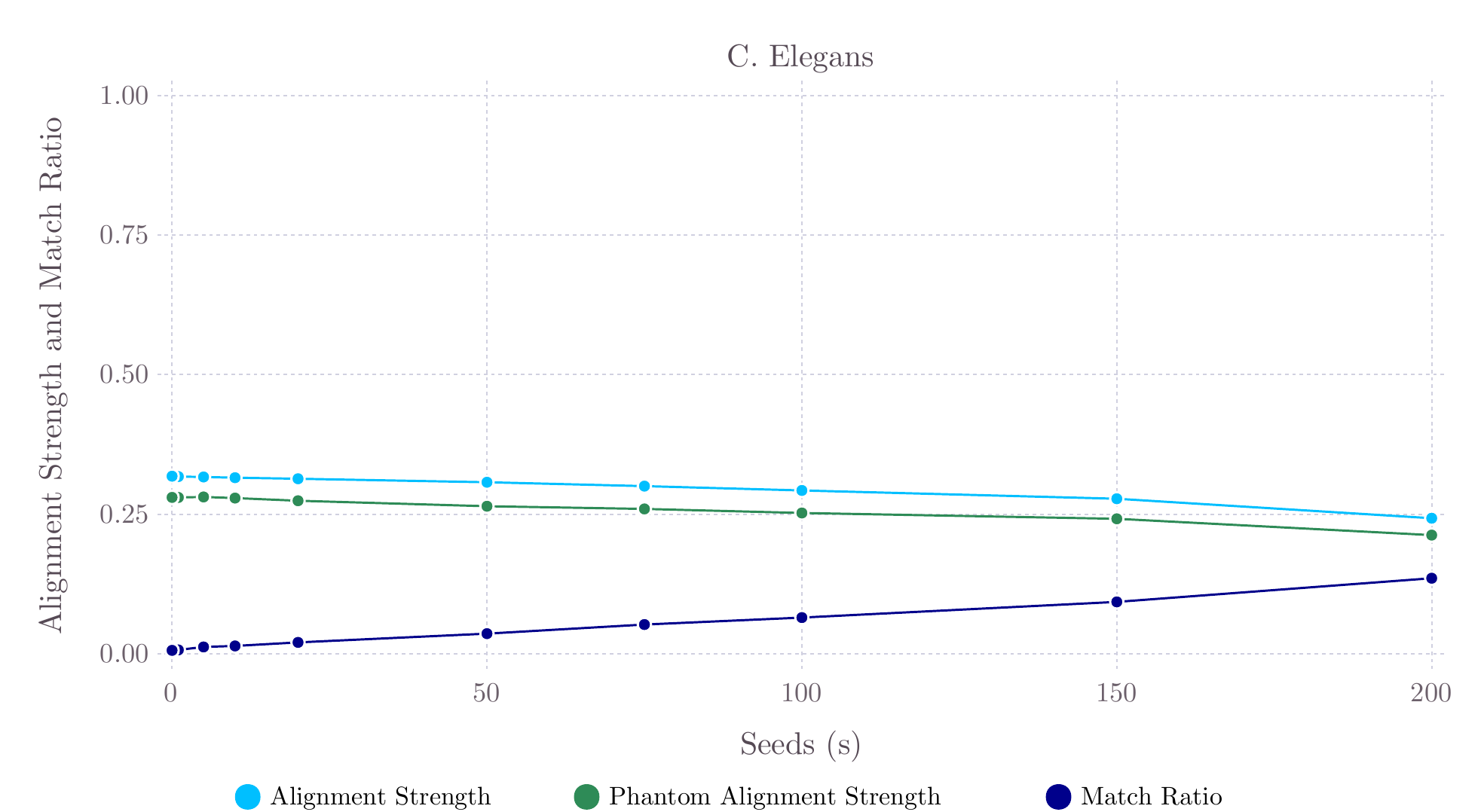}
    \caption{Matching C Elegans' electrical and chemical connectome}
    \label{fig:56:celegans}
\end{figure}

We conducted the identical experiments as we did for the
Wikipedia graphs, except that the number of seeds $s$
considered were $s=0, 1, 5, 10, 20, 50, 75, 100, 150, 200$, and we matched  $G^{e \ell}$ and  $G^{ch}$. The resulting plots are
in Figure \ref{fig:56:celegans};
alignment strength of the C Elegans
seeded graph match in blue,
phantom alignment strength in green, match
ratio in purple. Note that
seeded graph matching did poorly, as
evidenced by low match ratio, even when the number
of seeds was huge ($200$ seeds and $79$ nonseeds), and
correspondingly the gap between seeded graph matching alignment strength and phantom alignment strength was small.

The Enron graphs from  \cite{FAP} arose in the following
manner. Enron was a large and highly
respected energy company that dissolved
spectacularly in 2001 amid systemic fraud. The United
States Justice Department released a trove of email messages
between company employees. The graphs $G_{130}$, $G_{131}$,
and $G_{132}$ have as vertices $\n=184$ Enron employees and,
for each pair of vertices/employees,
the vertices are adjacent in $G_{130}$ precisely when
there is an email from one employee to the other in
week number 130 of the email corpus, they are adjacent
in $G_{131}$ precisely when there is an email from
one employee to the other in week number 131,
and they are adjacent
in $G_{132}$ precisely when there is an email from
one employee to the other in week number 132. The paper
\cite{Scanstat} identified an anomaly going into week 132, and
\cite{FAP} used match ratio differences between pairs
of these graphs to highlight this anomaly.

\begin{figure}[h!]
    \centering
    \includegraphics[width=0.6\textwidth]{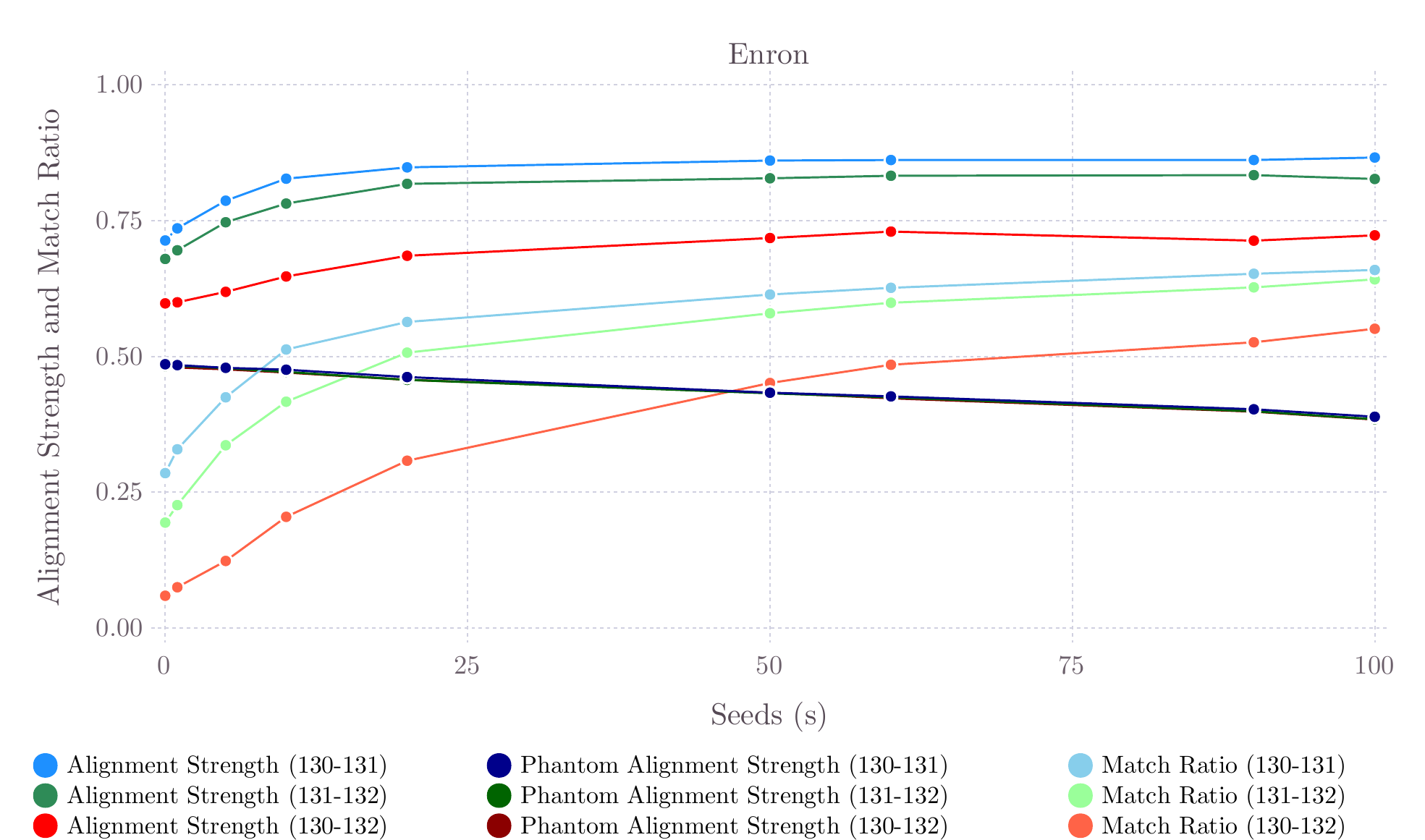}
    \caption{Matching Enron email networks; all pairs from weeks $\{130,131,132\}$.}
    \label{fig:56:enron}
\end{figure}

We conducted the identical experiments for each of the pairs
$G_{130},G_{131}$ and $G_{131},G_{132}$ and $G_{130},G_{132}$
as we did for the
Wikipedia graphs, except that  the number of seeds $s$
considered were $s=0, 1, 5, 10, 20, 50, 60, 90, 100$. The resulting plots are
in Figure \ref{fig:56:enron}. As noted in \cite{FAP},
the match ratio from matching $G_{130}$ to $G_{131}$ is
highest of the three, since the anomaly had not yet occurred.
The next highest match ratio was from matching $G_{131}$ to $G_{132}$, then came matching $G_{130}$ to $G_{132}$.
Note that the gap between seeded graph matching
alignment strength and phantom alignment strength was ordered
the same way; highest was $G_{130}$ to $G_{131}$, then
was $G_{131}$ to $G_{132}$, and then was $G_{130}$ to $G_{132}$. Indeed, more gap here when there was higher match ratio. (Note that the match ratios here differ a bit from those in the paper \cite{FAP}, Figure 8; that figure
was inadvertently from a nonsimple graph version of the data,
and here we created a simple graph.)

\section{Notable mentions and future directions, plus caveats \label{sec:notable} \label{sec:caveat}}

The applications of graph matching are broad and many, and
getting the right answer is only valuable when we know that
we have the right answer. This paper provides principled
tools that can help the practitioner decide if seeded graph
matching has found the true bijection.

The first caveat ---and future direction--- is that we are presenting a conjecture, and
not a theorem. Indeed, the Phantom Alignment Strength Conjecture, as formulated in Section \ref{sec:ASC}, includes
terms in quotes; ``moderate,'' ``high probability,'' ``very different,'' and
``negligible.'' Ironing these
terms out with specifics is part of the puzzle of
proving the conjecture, and is an important next task.
It may be a hard task, and we expect this paper to stimulate
more experimentation, fine-tuning, and eventually a proof
of the conjecture.

Part of the first caveat is the consideration
that we don't have a proof of the conjecture
as of now, and our experimentation is wide but not
exhaustive, and thus there may be additional hypotheses
or limitations to the conjecture statement.

 A second caveat
is that the conjecture is expressed in terms of an underlying
model for a pair of random graphs, and we need to
consider if particular real data that we may encounter (beyond the examples that we used here)
can more generally be considered as arising from such a model.

\begin{figure}[h!]
    \centering
    \includegraphics[width=0.7\textwidth]{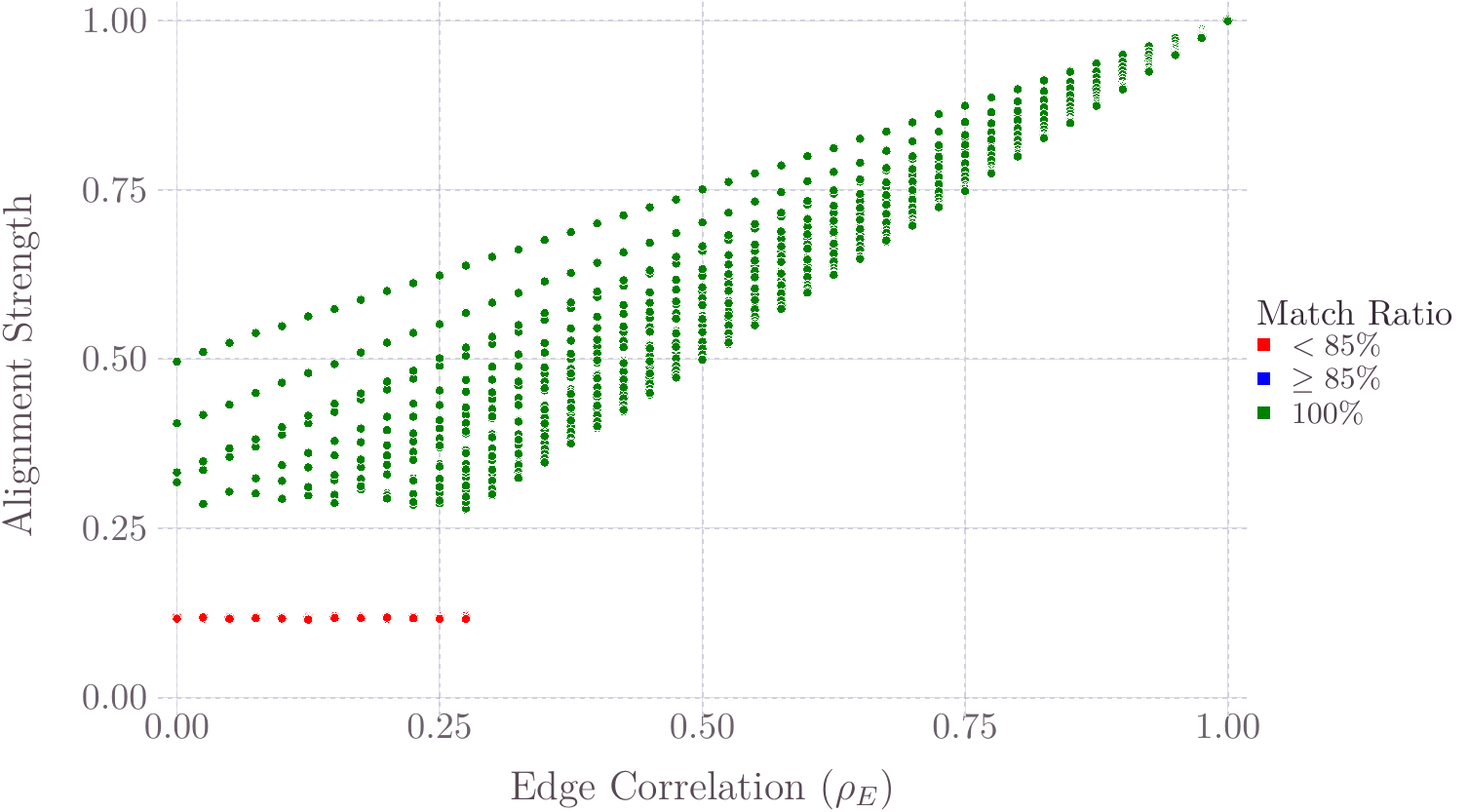}
    \caption{The same as the plot in  Figure \ref{fig:52a:d}, except that the x-axis in this figure is
    edge correlation instead of
    total correlation. The contrast between these two figures is quite stark, and highlights the utility of total correlation with regard to matchability. }
    \label{fig:52alt}
\end{figure}

Also, when there are multiple blocks with
Bernoulli coefficients being realized from different distributions for different blocks, we saw in Section \ref{sec:block}
that total correlation became a much less reliable
tool for determining matchability. More work is needed to
explore this further; the paper \cite{fishkind2019alignment},
when presenting empirical evidence
for the relationship between total correlation and
matchability, restricted their attention to the
setting hypothesized in our Phantom Alignment
Strength Conjecture, which excludes multiple blocks.
Indeed, in the setting of our conjecture, the role of
total correlation in matchability is starkly visible.
See Figure \ref{fig:52a:d}, where the x-axis is
total correlation, and compare to Figure \ref{fig:52alt}.
Figure \ref{fig:52alt} is the same data plotted in Figure \ref{fig:52a:d},
except that the x-axis is used
for edge correlation instead of total correlation.
The contrast between these two figures
is quite dramatic. Indeed, the current-literature-standard
yardstick of
edge correlation failed miserably in capturing
matchability, whereas total correlation captured matchability
perfectly here.
These two figures are powerful illustration of
the role of total correlation in matchability.

While more work remains to be done, we have here presented principled
tools that can be of significant help to the practitioner now.

A large measure of inspiration for this paper came from
Figures 1, 2, and 3 of  \cite{fishkind2019alignment} (on which half of us are co-authors). Those figures
displayed the results of graph matchings of many simulations of pairs of
correlated Bernoulli random graphs under similar conditions of the Phantom Alignment Strength Conjecture.
One axis of each figure tracked edge correlation, and the second axis tracked heterogeneity correlation; green, yellow, and red dots
were respectively located at coordinates corresponding to parameters where the graph matchings
were always the truth, mostly the truth, and often not the truth, respectively.
It was striking to observe that the regions of red and green were sharply
demarcated by a level curve of total correlation, with little yellow between the red and green.
These figures starkly demonstrated the role of total correlation
in matchability, as well as thresholding behavior. Together with the
theoretical results of \cite{fishkind2019alignment} tieing alignment
strength to total correlation (when graph matching gets truth),
we had important ingredients for the ``hockey stick'' at the heart of
the Phantom Alignment Strength Conjecture.
\\

\begin{backmatter}

\section*{Acknowledgements}
The authors thanks Dr. David Marchette for creating the Wikipedia graph
that we used here.

\section*{Funding}
This paper is based on research sponsored by the Air Force Research Laboratory and
DARPA under agreement number FA8750-20-2-1001 and FA8750-17-2-0112. The U.S. Government
is authorized to reproduce and distribute reprints for Governmental purposes notwithstanding any
copyright notation thereon. The views and conclusions contained herein are those of the authors
and should not be interpreted as necessarily representing the official policies or endorsements,
either expressed or implied, of the Air Force Research Laboratory and DARPA or the U.S. Government.

\section*{Abbreviations}
\noindent {\bf C. Elegans:} Caenorhabditis Elegans\\
\noindent {\bf dMRI:} Diffusion-weighted Magnetic Resonance Imaging\\
\noindent {\bf DT-MRI:} Diffusion Tensor Magnetic Resonance Imaging

\section*{Availability of data and materials}
Located at
\url{https://cs.jhu.edu/~fparker9/phantom-alignment-strength/}

\section*{Competing interests}
The authors declare that they have no competing interests.


\section*{Authors' contributions}
DEF formulated the Phantom Alignment Strength Conjecture, penned much of the manuscript.\\
FP, HS, and LM did the numerical experiments and contributed to the
formulation of the conjecture.\\
VL penned the introduction section, including the literature review, and made
the bibliography.\\
EB was involved in the acquisition of the human connectome data set, processed
the data for our particular use, and penned a portion of the manuscript involving the human and
C Elegans connectome networks.\\
VL, CEP, and AA worked on the theoretical development, edited the manuscript.\\
All authors read and approved the
manuscript.

\section*{Authors' information}


\bibliographystyle{bmc-mathphys} 
\bibliography{bibliot}      






\end{backmatter}
\end{document}